\documentclass{article}
\pdfoutput=1
\usepackage{amsthm}

\newtheorem{theorem}{Theorem}[section]
\newtheorem{lemma}[theorem]{Lemma}

\newtheorem{remark}[theorem]{Remark}

\usepackage{geometry}
\usepackage{graphicx}
\usepackage{epstopdf}
\usepackage{amscd}
\usepackage{amsmath,amssymb}
\usepackage{color}
\definecolor{cred}{rgb}{1,0,0}

  \newcommand{\p}{p}

  \newcommand{\pI}{p_I}
  \newcommand{\q}{q}

  \newcommand{\ph}{p_h}
   \newcommand{\qh}{\q_h}

  \newcommand{\fd}{\ell}
\newcommand{\z}{z}

 \renewcommand{\u}{{\boldsymbol u}}
  \newcommand{\bv}{{\boldsymbol v}}
  \newcommand{\w}{{\boldsymbol w}}
  \newcommand{\uh}{{\u}_h}
   \newcommand{\vh}{{\bv}_h}
  
  \renewcommand{\r}{r}
   \newcommand{\uI}{\u_{I}}
 
\newcommand{\bq}{{\boldsymbol{ q}}}
  \newcommand{\bw}{{\bf \w}}
  \newcommand\bg{\boldsymbol g}
\newcommand{\bchi}{{\boldsymbol \chi}}
\newcommand{\bzeta}{{\boldsymbol \zeta}}
\newcommand{\bphi}{{\boldsymbol \varphi}}
 \newcommand{\bui}{{\boldsymbol \Pi}^F_h{\u}}

  \newcommand{\UU}{{\bf U}}
  \newcommand{\WW}{{\bf W}}
  \newcommand{\VV}{{\bf V}}
  \newcommand{\FF}{{\bf F}}
  \newcommand{\ZZ}{{\bf Z}}
  \newcommand{\UUh}{{\bf U}_h}
  \newcommand{\VVh}{{\bf V}_h}
  \newcommand{\FFh}{{\bf F}_h}
  \newcommand{\bpsi}{\boldsymbol \Psi}

\newcommand{\bb}{{\bf b}}
\newcommand{\diff}{\nu}
\newcommand{\diffp}{{ \kappa}}

\newcommand{\vbeta}{\boldsymbol{\beta}}

\newcommand{\reaction}{\gamma}

\newcommand{\wiz}{\text{\tiny$\aleph$}}

  \newcommand{\n}{\boldsymbol{n}}
  \newcommand{\E}{E}
  \newcommand{\dE}{{\partial\E}}
  \newcommand{\Th}{{\mathcal T}_h}
  \newcommand{\Eh}{{\mathcal E}_h}
 \newcommand{\calO}{{\mathcal O}}
  \newcommand{\xx}{{\bf x}}
  \newcommand{\xb}{x_B}
  \newcommand{\yb}{y_B}
  \newcommand{\ds}{{\rm d}s}
  
 \newcommand{\dx}{{\rm d}x}

\newcommand{\RR}{ {\mathbb R}}

    \newcommand{\Pp}{\mathbb P}

    \newcommand{\VMhk}{{V}^k_h}
     \newcommand{\VMhkE}{{ V}^k_{h}(\E)}
    \newcommand{\V}{{\mathcal V}}
    \newcommand{\Vh}{{\mathcal V}_h}
    \newcommand{\calG}{ \mathcal{G}}

    \newcommand{\Amf}{{\mathfrak L}}
    \newcommand{\Amfs}{{\mathfrak L}^*}
    
    \renewcommand{\div}{\operatorname{div}}
    \newcommand\rot{\operatorname{rot}}

   \newcommand{\AAh}{{\mathcal A}_h}
\renewcommand{\AA}{{\mathcal A}}

      \newcommand{\NErrk}{{\mathcal E}^k}

    \newcommand{\Pzk}{{\Pi}^0_{k}}
    \newcommand{\PPzk}{{\bf\Pi}^0_{k}}

   \newcommand{\intE}{\int_\E}
   
   \newcommand{\br}[1]{\overline{#1}}

\begin{document}

\title{Mixed Virtual Element Methods for general second order elliptic problems on polygonal meshes}
%\thanks{...}\thanks{...}% At most 5 thanks
%
\author{L. Beir\~ao da Veiga\footnote{Dipartimento di Matematica, Universit\`a di Milano, Via
Saldini 50, 20133 Milano (Italy), and IMATI del CNR, Via Ferrata
1, 27100 Pavia, (Italy).},
%%\email{lourenco.beirao@unimi.it}
%%
\ F. Brezzi\footnote{IUSS, Piazza della Vittoria 15, 27100 Pavia
(Italy), and IMATI del CNR, Via Ferrata 1, 27100 Pavia, (Italy).},
%%\email{brezzi@imati.cnr.it}
%%
\ L.D. Marini\footnote{Dipartimento di Matematica, Universit\`a di
Pavia, and IMATI del CNR, Via Ferrata 1, 27100 Pavia, (Italy).},
%%\email{marini@imati.cnr.it}
%%
\ A. Russo\footnote{Dipartimento di Matematica e Applicazioni,
Universit\`a di Milano-Bicocca,via Cozzi 57, 20125 Milano (Italy)
and IMATI del CNR, Via Ferrata 1, 27100 Pavia, (Italy).}}
%%\email{alessandro.russo@unimib.it}

%
\date{...}

% \begin{resume} ... \end{resume}

%
\maketitle

%%%%%%%%%%%%%%%%%%%%%%%%%%%%%%%%%%%%%

%-%%%%%%%%%%%%%%%%%%%%--------------------------------------------------------------------
%-%%%%%%%%%%%%%%%%%%%%-
%-%%%%%%%%%%%%%%%%%%%%-
%-%%%%%%%%%%%%%%%%%%%%--------------------------------------------------------------------
%-%%%%%%%%%%%%%%%%%%%%-
%-%%%%%%%%%%%%%%%%%%%%-

{\bf Abstract:} In the present paper we introduce a Virtual
Element Method (VEM) for the approximate solution of general
linear second order elliptic problems in mixed form, allowing for
variable coefficients. We derive a theoretical convergence
analysis of the method and develop a set of numerical tests on a
benchmark problem with known solution.

% ----------------------------------------------------------------------------------------------------
\section{Introduction}\label{introduction}
% ----------------------------------------------------------------------------------------------------
%-%%%%%%%%%%%%%%%%%%%%--------------------------------------------------------------------
%-%%%%%%%%%%%%%%%%%%%%-
%-%%%%%%%%%%%%%%%%%%%%-
%-%%%%%%%%%%%%%%%%%%%%--------------------------------------------------------------------
%-%%%%%%%%%%%%%%%%%%%%-
%-%%%%%%%%%%%%%%%%%%%%-

The aim of this paper is to design and analyze some aspects of the use of  Virtual Element Methods (in short, VEM) for the approximate solution of general linear second order elliptic problems. In a previous paper \cite{variable-primal} the same authors analyzed diffusion-convection-reaction problems with variable coefficients in the primal form. Here we shall deal with the {\it mixed} formulation.

Virtual Element Methods (introduced in \cite{volley}) belong to the family of methods that allow the use of general polygonal and polyhedral decompositions, that are becoming more and more popular, in particular in view of their use in particular problems connected to moving boundaries.
Example of applications where polytopal meshes could have (or are already yielding) a positive impact can be found, for instance,
in fluid-structure interaction  \cite{XFEM84,WMLB},
crack propagation
\cite{Berrone-VEM,GFEM146,CRACKXFEM180},
phase change
\cite{Merle-Dolbow,chessa-S-B},
contact problems
\cite{Wriggers-2014},
or topology optimization
\cite{Gain-PhD,Paulino-VEM,TPPM10,Sutradhar-Paulino-Miller-Nguyen},
but they are promising also in other applications, for instance in presence of coefficients  that vary rapidly on sub-domains with complicated geometries, as when
dealing with various types of inclusions (see e.g. \cite{suku-CMB,Oswald,Paulino-nonlinear-polygonal}),
 or more generally in
medical applications \cite{SVC,vigneron,Oswald,spiegel},
in image processing
\cite{POLY28,HF,Floater-K-R},
and many others. It must be pointed out that several among these methods, in view of their great resistance to element distortions, come out to be handy not only for general polygonal elements, but also on quadrilaterals or hexahedra as well \cite{hourglass}.
The literature on these methods has quite old origins (see e.g.  \cite{Wachspress75}), and kept slowly increasing and widening its range of applications ever since. See for instance
\cite{ABCM,Arroyo-Ortiz,GFEM13,MESHLESS16,GFEM17,Bishop,Bochev-Hyman,DiPietro-Ern,Gillette-2,FHK06,Fries:Belytschko:2010,MESHLESS24,VEM19,GFEM-117,LIBROXFEM-128,Gillette-1,Suku04,Sukumar:Malsch:2006,ST04,POLY37,POLY47}.
In more recent times the variety of methods (already quite rich) has been growing very fast. In particular we have presently a flourishing group of methods, quite similar to each other, based (one way or another) on local polynomial reconstructions. Among others we mention Hybridizable Discontinuous Galerkin methods
\cite{Cockburn-IMU,Cockburn-Jay-Lazarov,Cockburn-Jay-Sayas,Cockburn-Guzman-HWang,Nguyen-Peraire-Cockburn},
Weak Galerkin methods
\cite{Mu-Wang-Wei-Ye-Zhao,Mu-Wang-Ye-2,Wang-1,Wang-2},
the latest evolution of Mimetic Finite Differences
\cite{BLM09,BLM11,BeiraodaVeiga:Manzini:2008b,Brezzi:Buffa:Lipnikov:2009,VEM19},
several variants related to Finite Volumes and Mixed
Methods \cite{Droniou-SpIss,Droniou-gradient,Droniou:Eymard:Gallouet:Herbin:2010,Bonelle:Ern:2014,DiPietro-Ern-1,DiPietro-Ern-2,DiPietro-Ern-3,DiPietro-Ern--Lemaire},
boundary element methods \cite{BEM-weisser}
and various evolutions of the Virtual Element Methods themselves (mentioned below).

The similarities and the differences among all these methods are still under investigation, as well as the (much more important) analysis of  ``which method is best suited for which class of problems". We are not going to attempt to clarify these issues in the present paper, and more modestly we stick on Virtual Element Methods, and in particular on their use in mixed formulations.

We recall that Mixed Virtual Element Methods for $\div({\mathbb K}\nabla)$ with ${\mathbb K}$ constant
were introduced, for the two dimensional case, in \cite{BFM}  as an evolution of the Mimetic Finite Differences as originally analyzed in \cite{Brezzi:Lipnikov:Shashkov:2005,Brezzi:Lipnikov:Shashkov:Simoncini:2007,Brezzi:Lipnikov:Shashkov:2006,Brezzi:Lipnikov:Simoncini:2005}, and then extended in various directions, see for instance \cite{Beirao-apos,Maxwell-MFD,MFD7,Antonietti}.
For references to several much older papers on Mimetic Finite Differences and a much more detailed panorama  on related methods  we refer  to \cite{MFD22} and \cite{MFD23}. We also point out that the first attempt to extend and analyze Mimetic Finite Differences to linear elliptic  second order operators of the form $\div({\mathbb K}\nabla)$  with a variable ${\mathbb K}$ was actually  done earlier in \cite{BLM09} for the mixed formulation.

A more recent approach to the theory of Virtual Element Methods has been introduced  in \cite{projectors}, where the first
attempt to a systematic use of the $L^2$-projection operator was presented (originally for the so called {\it nodal} VEM). This was later refined and extended to mixed formulations in \cite{super-misti}. See also
\cite{hitchhikers},
for more details on the implementation of Virtual Elements and \cite{Brezzi:Marini:plates,VEM-elasticity,BM13,Berrone-VEM,Manzini:Russo:Sukumar,Paulino-VEM,ABMV14,VemSteklov} for other  interesting applications and developments.

Here we follow this direction, and the Virtual Element Methods that we propose and analyze for dealing with variable coefficients are indeed
based on $L^2$-projection operators in a rather systematic way. We recall that for Virtual Element Methods the shape and trial functions
are not given in an explicit form, but rather as solutions of  PDE problems inside each element. As we do not want to solve these problems
inside the elements (not even in an approximate way), the passage from constant to variable coefficients is less trivial than for other methods. In particular, simple minded approaches to variable coefficients can lead to a loss of optimality, especially for higher order methods, as it has been shown for instance in  \cite{variable-primal} for nodal VEM.

For the sake of simplicity we present here only the two-dimensional case, although, as pointed out here below in Remark
\ref{2-3Dim}, the passage from two to three dimensions, in the present case, is quite immediate.

We will use the following notation. The space of polynomials of degree $\le k$, for $k$ nonnegative integer, will be denoted by $\Pp_k$, or $\Pp_k(\calO)$ whenever we want to stress the fact that we are working on a particular domain $\calO$. As common, we will use $\Pp_{-1}\equiv \{0\}$ as well.

Throughout the paper, we will follow the standard notation for classical Sobolev spaces, as for instance in \cite{Ciarlet-78}. In particular, for a domain $\calO$ in one or several dimensions,
$\| f\|_{k,p,\calO}$ ($k\ge 0$ integer and $1\le\p\le+\infty$) will denote the norm of the function $f$ in the Sobolev space $W^{k,p}(\calO)$ of functions that belong to $L^p(\calO)$ with all their derivatives up to the order $k$. We will also use the notation
$H^k(\calO)$ to denote $W^{k,2}(\calO)$, and the norm of a function $f$ in $H^k(\calO)$ will be denoted by
$\|f\|_{k,\calO}$ (or simply $\|f\|_k$ whenever no confusion can occur). With a minor (and common) abuse of notation, for a vector valued function (say, ${\bf f}:\calO\rightarrow \RR^2$) we will still write  $\|{\bf f}\|_{k,p,\calO}$ to denote the norm of  ${\bf f}$ in the Sobolev space $(W^{k,p}(\calO))^2$.
%For a domain $\calO$ (in one or two dimensions),
The scalar product in $L^2(\calO)$ or in $(L^2(\calO))^2$ will be denoted  by   $(\cdot\,,\,\cdot)_{0,\calO}$, or simply by   $(\cdot\,,\,\cdot)_0$ (or even
 $(\cdot\,,\,\cdot)$) when no confusion may arise. As usual, $H^k_0(\calO)$ ($k$ integer $>0$) will denote the subset of $H^k(\calO)$
made of functions vanishing at the boundary $\partial\calO$ of $\calO$ together with all their derivatives up to the order $k-1$.

 Throughout the paper, $C$ will denote a generic constant independent of the mesh size, not necessarily the same from one occurrence to the other. Sometimes, in some specific step where we want  to stress the dependence of a constant on some variable (say, $\xi$) we will indicate it by $C_{\xi}$. Needless to say, $C_{\xi}$ might also assume different values from one occurrence to another.

An outline of the paper is as follows.  In Section \ref{problem}, after stating the problem and its formal adjoint, we recall (in Subsection \ref{sec:mixed})  the mixed variational formulation. Then, in Section \ref{VEM-mixed} we introduce the Virtual Element approximation of the mixed formulation, and derive optimal error estimates in Section \ref{sec:mixed:est}.  In Section \ref{superconvergence} we derive a superconvergence result for the scalar variable, and finally, in Section \ref{Num-Exp}, we present some numerical results.

In the bibliography we included an unusual amount of references, as it would have been appropriate for a review paper. However we thought that a wide set of references could be convenient, as well,  for a paper submitted {for} a special issue (like the present one).

%-%%%%%%%%%%%%%%%%%%%%--------------------------------------------------------------------
%-%%%%%%%%%%%%%%%%%%%%-
%-%%%%%%%%%%%%%%%%%%%%-
%-%%%%%%%%%%%%%%%%%%%%--------------------------------------------------------------------
%-%%%%%%%%%%%%%%%%%%%%-
%-%%%%%%%%%%%%%%%%%%%%-
% ----------------------------------------------------------------------------------------------------
\section{The problem and the adjoint problem}\label{problem}
% ----------------------------------------------------------------------------------------------------
%-%%%%%%%%%%%%%%%%%%%%--------------------------------------------------------------------
%-%%%%%%%%%%%%%%%%%%%%-
%-%%%%%%%%%%%%%%%%%%%%-
%-%%%%%%%%%%%%%%%%%%%%--------------------------------------------------------------------
%-%%%%%%%%%%%%%%%%%%%%-
%-%%%%%%%%%%%%%%%%%%%%-

Let $\Omega\subset \RR^2$ be a bounded convex polygonal domain and let $\Gamma$ represent the boundary of $\Omega$. We assume that  $\diffp$  and $\reaction$ are smooth functions $\Omega\rightarrow\RR$  with $\diffp(\xx)\ge \diffp_0>0$ for all $\xx\in\Omega$, and that $\bb$  is a smooth vector valued function $\Omega\rightarrow\RR^2$.
For $f\in H^{-1}(\Omega) (\equiv (H^1_0(\Omega))^\prime)$, we consider the problem:
\begin{equation}\label{Pb-cont}
\left\{
\begin{aligned}
&\mbox{Find $\p\in H^1_0(\Omega)$ such that:}\\
&\Amf\,\p:=\div(-\diffp(\xx)\nabla\p + \bb(\xx)\p) + \reaction(\xx)\,\p
= f(\xx)\quad\text{\it  in }\Omega.\\
\end{aligned}
\right.
\end{equation}
We make the following fundamental assumption, that among other things implies that problem \eqref{Pb-cont} is Well-Posed.

\noindent {\bf Assumption WP} We assume that for all source terms $f\in H^{-1}(\Omega)$  problem \eqref{Pb-cont}
has a unique solution $\p$, that moreover satisfies the a-priori estimate
\begin{equation}\label{bound11}
\|\p\|_{1,\Omega} \le C \|f\|_{-1,\Omega} ,
\end{equation}
as well as the regularity estimate
\begin{equation}\label{bound02}
\|\p\|_{2,\Omega} \le C \|f\|_{0,\Omega},
\end{equation}
both with a constant $C$ independent of $f$.

\medskip
We consider also the adjoint operator $\Amfs$ given by
\begin{equation}\label{def:adjoint}
\Amfs\p:=\div(-\diffp(\xx)\nabla\p)- \bb(\xx)\cdot\nabla\p + \reaction(\xx)\,\p.
\end{equation}
The above assumptions on problem \eqref{Pb-cont} imply, among other things, that
existence and uniqueness hold, as well, for \eqref{def:adjoint}.
Moreover, for every $g\in L^2(\Omega)$ there exists a unique $\varphi\in H^2(\Omega)\cap H^1_0(\Omega)$ such that $\Amfs\varphi=g$, and
\begin{equation}\label{stimadj}
\|\varphi\|_{2,\Omega} \le C^* \|g\|_{0,\Omega}
\end{equation}
for a constant $C^*$ independent of $g$.
We note that having a full diffusion tensor would not change the analysis in a substantial way; the choice of having a scalar diffusion coefficient $\diffp$ was done just for simplicity.
Finally, as we shall see, the 2-regularity \eqref{bound02} and \eqref{stimadj} is not necessary in order to
derive the results of the present work, and an $s$-regularity with $s>1$ would be sufficient. Here however we are not interested in minimizing the regularity assumptions.

% Actually, in the numerical results presented in Section \ref{Num-Exp} a full tensor is used.

%%%%%%%%%%%%%%%%%%%%%%%%%%%
%%%%%%%%%%%%%%%%%%%%%%%%%%%
\subsection{The mixed variational formulation}\label{sec:mixed}
%%%%%%%%%%%%%%%%%%%%%%%%%%%
%%%%%%%%%%%%%%%%%%%%%%%%%%%

In order to build the mixed variational formulation of problem \eqref{Pb-cont}, we define
$$ \diff:=\diffp^{-1},\quad \vbeta:=\diffp^{-1} \bb,$$
and re-write \eqref{Pb-cont} as
\begin{equation}\label{splitted}
 \u= \diff^{-1}(-\nabla\p + \vbeta \p), \quad \div \u +\reaction \, \p =f~~\text{in } \Omega,\quad \p=0~~\text{on }\Gamma.
\end{equation}
Introducing the spaces
$$ V:= H(\div;\Omega), \quad \text{and}\quad Q:= L^2(\Omega),$$
the variational formulation of problem \eqref{splitted} is:
\begin{equation}\label{mixed-cont}
\left\{
\begin{aligned}
&\text{Find~} (\u,\p)\in V\times Q\text{ such that}\\
&(\diff \u, \bv)-(\p,\div \bv)-(\vbeta\cdot \bv, \p)= 0\quad \forall \bv \in V,\\
&(\div \u, q)+(\reaction \p,q)=(f,q)\quad \forall q\in Q.
\end{aligned}
\right.
\end{equation}
For the subsequent analysis it will be convenient to write \eqref{mixed-cont} also  in a more compact way.
For this, we define first
$$\V:=V\times Q, \quad
%and
%$$
\UU:=(\u,\p),\quad \VV:=(\bv,q),\quad \FF:=(0,f),
$$
and
\begin{equation}\label{def:AA0}
\AA(\UU,\VV):=(\diff \u, \bv)-(\p,\div \bv)-(\vbeta\cdot \bv, \p)
+(\div \u, q)+(\reaction \p,q).
\end{equation}
Problem \eqref{mixed-cont} can then be equivalently written as:
\begin{equation}\label{def:AA}
\left\{
\begin{aligned}
&\mbox{Find } \UU\in \V \mbox{ such that}\\
&\AA(\UU,\VV)=(\FF,\VV)\quad \forall \VV\in \V.
\end{aligned}
\right.
\end{equation}

\begin{remark} It is almost immediate to see that our path (from \eqref{Pb-cont}) to \eqref{def:AA}) can be easily reversed:
if a pair $\UU=(\u,\p)$ solves \eqref{def:AA} then $\u$ and $\p$ satisfy \eqref{splitted} and hence $\p$ solves \eqref{Pb-cont}.
In turn, this easily gives that the existence and uniqueness of the solution of \eqref{Pb-cont} implies the existence and uniqueness of the solution of \eqref{def:AA}.
\end{remark}

%-%%%%%%%%%%%%%%%%%%%%--------------------------------------------------------------------
%-%%%%%%%%%%%%%%%%%%%%-
%-%%%%%%%%%%%%%%%%%%%%-
\section{VEM approximation}\label{VEM-mixed}
%-%%%%%%%%%%%%%%%%%%%%-
%-%%%%%%%%%%%%%%%%%%%%-
%-%%%%%%%%%%%%%%%%%%%%-

In the present section we introduce the Virtual Element approximation of problem \eqref{mixed-cont}.

%-%%%%%%%%%%%%%%%%%%%%--------------------------------------------------------------------
%-%%%%%%%%%%%%%%%%%%%%-
%-%%%%%%%%%%%%%%%%%%%%-
\subsection{The Virtual Element spaces}

%-%%%%%%%%%%%%%%%%%%%%--------------------------------------------------------------------
%-%%%%%%%%%%%%%%%%%%%%-
%-%%%%%%%%%%%%%%%%%%%%-

Let $\Th$ be a decomposition of $\Omega$ into star-shaped polygons $\E$, and let $\Eh$ be the set of edges $e$ of $\Th$.
We {further} assume that {for every element $\E$  there exists a $\rho^{\E}>0$ such that $\E$}  is star-shaped with respect to every point of a disk $D_{{\rho_{\E}}}$ of radius
$\rho^{\E} h_E$ (where $h_{\E}$ is the diameter of $\E$) and that the length $h_e$ of every edge $e$ of $\E$ satisfies $h_e \ge  {\rho^{\E}} h_{\E}$.
When considering a {\it sequence} of decompositions $\{\Th\}_h$ we will obviously assume $\rho^{\E}\ge\rho_0>0$ for some $\rho_0$ independent of $\E$ and of the decomposition. As usual, $h$ will denote the maximum diameter of the elements of $\Th$.

For every element $\E$ we introduce:
\begin{equation} \label{defG}
\calG_{k}(\E):=\nabla \Pp_{k+1}(\E),
\end{equation}
and
%we denote
\begin{equation}\label{defGperp}
\calG^\perp_{k}(\E)=\text{  the $L^2(\E)$ orthogonal of  } \calG_{k}(\E) \text{ in }(\Pp_k(\E))^2,
\end{equation}
so that
\begin{equation}\label{oplus}
(\Pp_{k}(\E))^2=\calG_{k}(\E) \oplus \calG^\perp_{k}(\E).
\end{equation}
%
 %On each element $\E\in\Th$
 For $k$ integer $\ge 0$ we define
\begin{multline}\label{spazio-VE}
\VMhkE:=\{\bv\in H(\div;\E)\cap H(\rot;\E):
 \bv\cdot\n_{|e}
\in\Pp_{k}(e)~\forall  e
 \in \dE,\\
\div\bv\in\Pp_{k}(\E), \mbox{ and } \rot\bv\in\Pp_{k-1}(E)\}.
\end{multline}
Then we introduce the discrete spaces
\begin{equation}\label{spazio-V}
{\VMhk}:=\{\bv\in H(\div;\Omega) \mbox{ such that } \bv_{|\E}\in \VMhkE ~\forall \mbox{ element $\E$ in } \Th\},
\end{equation}
and
\begin{equation}\label{spazio-Q}
Q^k_h:=\{q\in L^2(\Omega) \mbox{ such that: } q_{|\E}\in \Pp_{k}(\E)~\forall \mbox{ element $\E$ in }
\Th\}.
\end{equation}
%{\color{red}
%\begin{remark}\label{reg-trace} We point out that the regularity  $\bv\in H(\div;\E)\cap H(\rot;\E)$  required in \eqref{spazio-VE}
%is enough to allow us consider the restriction of $\bv\cdot\n_{|e}$ to every single edge $e$  (see e.g. \cite{ costabel} ) while
%the regularity $\bv\in H(\div;\E)$, as is well known, would not have been enough for that.
%\end{remark}
%}
The degrees of freedom for $Q^k_h$ are obvious (one has many equivalent good choices for them),  while the degrees of freedom for $\VMhk$ are defined by (see \cite{super-misti})
\begin{eqnarray}%\label{doffdd0}
&\int_e{\bv\cdot\n}\,{q}_{\,k}\,\ds
&\quad \mbox{ for all edge $e$, for all }\;
q_{k}\in\Pp_k(e),\label{dof1}\\[3pt]
&\int_{\E}{\bv\cdot \bg_{k-1}} \dx &\quad \mbox{ for all element $\E$, for all $\bg_{k-1}\in\calG_{k-1}(\E)$}, \quad \label{dof2}
\\[3pt]
&\int_{\E}{\bv \cdot\bg^\perp_{k}} \dx  &\quad \mbox{ for all element $\E$, for all $\bg^\perp_{k}\in \calG^\perp_{k}(\E)$},
\label{dof3}
\end{eqnarray}
where the notation \eqref{defG}-\eqref{defGperp} was used  for  $\calG_{k}(\E)$  and   $\calG^\perp_{k}(\E)$,  respectively.
\begin{remark}
We point out that conditions \eqref{dof1} could be replaced by the values of $\bv \cdot \n$ at suitable points on each edge. Similarly, in \eqref{dof3} $\calG^\perp_{k}(\E)$ could be replaced by any subspace of $(\Pp_k(\E))^2$ satisfying \eqref{oplus}.
\end{remark}
\begin{remark} It is not difficult to check that the present choice of elements mimics, in some sense,  the Raviart-Thomas elements,
although, even on triangles, they coincide with the RT elements only for $k=0$. As pointed out in \cite{BFM} and in \cite{super-misti}
there are many other choices that could be made.
\end{remark}

\begin{remark}
Regarding the mesh assumptions at the beginning of this section, we note that it wouldn't be a problem to generalize the shape regularity condition by allowing suitable unions of star-shaped elements. Analogously, also the minimal edge length condition could be probably avoided with some additional technical work in the interpolation estimates.
\end{remark}

%-%%%%%%%%%%%%%%%%%%%%--------------------------------------------------------------------
%-%%%%%%%%%%%%%%%%%%%%-
%-%%%%%%%%%%%%%%%%%%%%-
\subsection{Interpolants, projections and approximation errors}
%-%%%%%%%%%%%%%%%%%%%%--------------------------------------------------------------------
%-%%%%%%%%%%%%%%%%%%%%-
%-%%%%%%%%%%%%%%%%%%%%-

From now on, we shall denote by $\Pzk:~Q\rightarrow Q^k_h$ and by $\PPzk:~V\rightarrow {\VMhk}$ the $L^2-$ projection
operators, defined locally by
\begin{equation}\label{projections-2}
\begin{aligned}
&\int_{\E}(q-\Pzk q) p_k\,dx =0~\forall p_k \in \Pp_k(\E),\quad \forall \E\in \Th,\\
&\int_{\E}(\bv-\PPzk \bv){\bq}_{k}\,dx =0~\forall {\bq}_{k} \in (\Pp_k(\E))^2 ,\quad \forall \E\in \Th.
\end{aligned}
\end{equation}

In \cite{super-misti} it was shown that the degrees of freedom \eqref{dof1}-\eqref{dof3} allow the explicit computation of the
projection $\PPzk \bv$ from the knowledge of the degrees of freedom \eqref{dof1}-\eqref{dof3} of $\bv$.  For the convenience of the reader we briefly recall the construction. We first observe that
using the degrees of freedom \eqref{dof2} we can easily compute the value of $\bv\cdot\n$ on $\partial\E$. From this and
\eqref{dof2} one can compute the value of $\div\bv\in\Pp_k$, using
\begin{equation}\label{calcdiv}
\int_{\E}{\div\bv\, q_k} \dx=-\int_{\E}{\bv\cdot \nabla q_k}\dx+\int_{\partial\E}{\bv\cdot \n \,q_k}\ds\qquad\forall q_k\in\Pp_k
\end{equation}
(remember that $\nabla q_k\in\calG_{k-1}$).
Once you know explicitly $\bv\cdot\n$ on $\partial\E$ and $\div\bv$ inside $\E$, then you can easily compute the integral
\begin{equation}
\int_{\E}{\bv\cdot \nabla q_{k+1}}\dx=-\int_{\E}{\div\bv \,q_{k+1}} \dx+\int_{\E}{\bv\cdot \n\, q_{k+1}}\ds,
\end{equation}
meaning that you can compute $\int_{\E}{\bv\cdot \bg_{k}} \dx$ for every $\bg_k\in\calG_k$. This and the degrees of freedom
\eqref{dof3} allow you to compute  $\int_{\E}{\bv\cdot \bq_{k}} \dx$ for every (vector valued) polynomial $\bq_k$ of degree $\le k$.

On the other hand, in every element $\E$, the computation of the $L^2(\E)$-projection of an element $q\in Q^k_h$ is trivial
(and coincides with its restriction to the element $\E$).

With classical arguments one can easily show that
\begin{equation}\label{projection-error}
\|q-\Pzk q\|_0\le C h^{s} |q|_{s},\quad \| \bv-\PPzk \bv\|_0 \le C h^{s} |\bv|_{s},~0\le s\le k+1,
\end{equation}
for every $q$ and $\bv$, respectively, that make the norms in the right-hand sides finite.

We point out that a linear ``Fortin'' operator  ${\boldsymbol\Pi}^F_h$ from $W:=(H^1(\Omega))^2\rightarrow \VMhk$ can be defined through the degrees of freedom \eqref{dof1}-\eqref{dof3}, by setting, brutally
\begin{eqnarray}%\label{doffdd0}
&\int_e{(\bv-{\boldsymbol\Pi}^F_h\bv)\cdot\n}\,{q}_{\,k}\,\ds=0
&\quad \mbox{ for all edge $e$, for all }\;
q_{k}\in\Pp_k(e),\label{doff1F}\\[3pt]
&\int_{\E}{(\bv-{\boldsymbol\Pi}^F_h\bv)\cdot \bg_{k-1}} \dx=0 &\quad \mbox{ for all element $\E$, for all $\bg_{k-1}\in\calG_{k-1}(\E)$}, \quad \label{doff2}
\\[3pt]
&\int_{\E}{(\bv-{\boldsymbol\Pi}^F_h\bv) \cdot\bg^\perp_{k}} \dx =0 &\quad \mbox{ for all element $\E$, for all $\bg^\perp_{k}\in \calG^\perp_{k}(\E)$},
\label{doff3}
\end{eqnarray}
and (using, essentially, \eqref{calcdiv}) it is easy to verify that the commuting diagram property holds:
\begin{equation}\label{commuting}
\begin{CD}
W                          @>{\displaystyle\div}>>       Q        @>>>0                \\
@V{\displaystyle{\boldsymbol\Pi}^F_h}VV                          @VV{\displaystyle\Pzk}V       \\
{\displaystyle\VMhk}       @>>{\displaystyle\div}>       Q^k_h               @>>>0
\end{CD}
\end{equation}
so that
\begin{equation}\label{prop-div}
\div {\boldsymbol\Pi}^F_h \bv = \Pzk \div \bv.
\end{equation}
Moreover, the following estimates hold, provided  $\u$ has enough regularity:
\begin{equation}\label{interp-VEM}
%\begin{aligned}
\|\u-\bui\|_0 \le C h^{k+1}\|\u\|_{k+1},\quad
%\|\p-\Pzk\p\|_0\le C h^{k+1}\|\p\|_{k+1},\\
\|\div(\u-\bui)\|_0 \le C h^{k+1}\|\div\u\|_{k+1}    .
%\end{aligned}
\end{equation}
With a minor abuse of notation, for an element $\WW\equiv (\w,r)$ with $ \w \in (H^1_0(\Omega))^2$ and $r$ scalar or vector function in  $L^2(\Omega)$, we will also denote
\begin{equation}\label{medie}
\begin{aligned}
&\overline{\w}:=\PPzk\w,\qquad \overline{r}:=\Pzk r,\quad\mbox{and}\quad
\overline{\WW}:=(\overline{w},\overline{r}),\\
&\w_I:=\Pi^F_h \w,\qquad r_I:=\Pzk r,\quad\mbox{and}\quad \WW_{I}:=(\w_I,r_I).
\end{aligned}
\end{equation}
We remind that, obviously,
\begin{equation}\label{projred}
\|\overline{\w}\|_0\le\,\|\w\|_0,\quad\quad\|\overline{r}\|_0\le\,\|r\|_0,\quad\quad\|r_I\|_0\le\,\|r\|_0,
\end{equation}
while
\begin{equation}\label{bound_I}
\|\w_I\|_0\le\,\|\w\|_0 +\|\w-\w_I\|_0\le C\,(\|\w\|_0+ h |\w|_1).
\end{equation}

For locally smooth $\w$, as we can see from \eqref{projection-error} and \eqref{interp-VEM}, the two errors
$\|\w-\w_I\|_{0,\E}$ and $\|\w-\br{\w}\|_{0,\E}$ will behave in the same way (in terms of
powers of $h$ and required regularity). Hence it makes sense to introduce a sort of {\it common value} that bounds both
of them. We define
\begin{equation}\label{defErrw}
\NErrk(\w):=\|\w-\w_I\|_{0}+\|\w-\br{\w}\|_{0},
\end{equation}
and we put it on charge to measure the {\it approximation error} for $\w$ when using Virtual Element spaces
of degree $k$. Needless to say, the same holds for a scalar function $r$ approximated in $Q^k_h $ (or, when
necessary, for a vector valued function ${\bf r}$ approximated in $(Q^k_h )^2$)  since, in these cases
the two approximations $\br{r}$ and $r_I$ coincide (as one can see in \eqref{medie}). In order to use
the same notation all over, however, we follow \eqref{defErrw}, and set
\begin{equation}\label{defErrrW}
\NErrk(r):=\|r-r_I\|_{0}+\|r-\br{r}\|_{0},\quad \mbox{and}\quad \NErrk(\WW):=\|\WW-\WW_I\|_{0}+\|\WW-\br{\WW}\|_{0}.
\end{equation}
We also point out that, by the properties of the projection, we immediately have
\begin{equation}\label{projWI}
\|\WW_I-\br{\WW_I}\|_{0}\le\|\WW_I-\br{\WW}\|_{0}\le \|\WW_I-{\WW}\|_{0}+\|\WW-\br{\WW}\|_{0}
= \NErrk(\WW),
\end{equation}
implying also
\begin{equation} \label{errWI}
\NErrk(\WW_I)\le \NErrk(\WW).
\end{equation}

%{\color{cred}
%   \begin{remark} It is obvious that the second term in \eqref{defErrw} can be bounded in %terms of the first (with constant equal to one) so that all the discussion concerning the %symbol $\NErrk$, so far, looks quite useless. This however will not be the case when a %variable coefficient  is also taken into account, as in the following notation (for instance %in \eqref{Ewrfi}-\eqref{selava}) that will help allowing a shorter treatment of the error %estimates in the next sections
%\qed
%\end{remark}
%}

Along the same lines, it is intuitively obvious (and it can be easily proved) that if you have a certain estimate  (in terms of powers of $h$ and required regularity) for $\w$ (or for $r$) you will have  quite similar estimates for, say,
$\varphi\w$ whenever $\varphi$ is a given smooth function. The constant in front of the estimate will
depend on $\varphi$, but the power of $h$ and the regularity required to $\w$ will be exactly the same. For instance it is immediate to check (just expanding the derivatives of the products, and using Cauchy-Schwarz) that one has
\begin{equation}
\|\varphi\w-\br{\varphi\w}\|_0\le C\,h^{k+1}\,|\varphi\w|_{k+1}\le C\,h^{k+1}\,\|\varphi\|_{k+1,\infty}\|\w\|_{k+1}
\equiv C_{\varphi}\,h^{k+1}\,\|\w\|_{k+1}.
\end{equation}

The same occurs for a pair $\WW=(\w,r)$  when one of the two entries (or both) are multipled by a smooth
function $\varphi$ or a smooth vector valued function ${\bphi}$, as in
\begin{equation}\label{Ewrfi}
\NErrk(\w\varphi)=\|\w\varphi-\br{\w\varphi}\|_0+\|\w\varphi-({\w\varphi})_I\|_0
\quad\mbox{and}\quad
\NErrk(r\varphi)=\|r\varphi-\br{r\varphi}\|_0+\|r\varphi-({r\varphi})_I\|_0
\end{equation}
for a smooth function $\varphi$, as well as in
\begin{equation}\label{Erbfi}
\NErrk(r\bphi)=\|r\bphi-\br{r\bphi}\|_0+\|r\bphi-({r\bphi})_I\|_0
\end{equation}
for a smooth vector valued function $\bphi$.

All this suggests a further ``abuse of notation'':  for  $\WW=(\w,\r)$ we will use the notation  $\NErrk(\wiz\WW)$
(either for $\wiz$ scalar or $\wiz$ vector)  whenever one of the two ($\w$ and  $r$), or both, are multiplied by
$\wiz$.

It could be worth pointing out a few particular cases: no matter whether $\wiz$  is a scalar or  a vector, we have
\begin{equation}\label{inH1}
\NErrk(\wiz\WW)\le C_{\wiz}\,h\,\Big(\|\r\|_1+\|\w\|_1\Big),
\end{equation}
as well as
\begin{equation}\label{selava}
\NErrk(\wiz\WW_I)\le\NErrk(\wiz(\WW_I-\WW))+\NErrk(\wiz\WW)
\le \|\wiz\|_{\infty}\NErrk(\WW)+\NErrk(\wiz\WW).
\end{equation}
Needless to say, the obvious analog of the bounds \eqref{inH1}-\eqref{selava} apply also to the separate terms $\NErrk(\w),~\NErrk(\r)$
%\begin{equation}
%\NErrk(\w):=\|\w-\w_I\|_0+\|\w-\br{\w}\|_0\quad\mbox{and}\quad\NErrk(\r):=\|\r-\br{\r}\|_0,
%\end{equation}
and so on. Finally, we observe that estimates \eqref{projection-error} and \eqref{interp-VEM} imply
\begin{equation}\label{stime-variE}
\NErrk(\UU)\le C h^{k+1}(\|\u\|_{k+1}+\|\p\|_{k+1}),\quad \NErrk(\wiz\UU)\le C_{\wiz}h^{k+1}(\|\u\|_{k+1}+\|\p\|_{k+1}),
\end{equation}
where $C_{\wiz}$ is a constant depending on $\wiz$ and its derivatives up to the order $k+1$.
As a final remark we note that, {\it whenever convenient},
we can easily bound $\NErrk(\wiz\WW)$ by
\begin{equation}
\label{consi6}
\NErrk(\wiz\WW)\le C_{\wiz}\|\WW\|_0 .
\end{equation}

%-%%%%%%%%%%%%%%%%%%%%--------------------------------------------------------------------
%-%%%%%%%%%%%%%%%%%%%%-
%-%%%%%%%%%%%%%%%%%%%%-
\subsection{The discrete bilinear forms}
%-%%%%%%%%%%%%%%%%%%%%--------------------------------------------------------------------
%-%%%%%%%%%%%%%%%%%%%%-
%-%%%%%%%%%%%%%%%%%%%%-

As is well known from the theory of mixed formulations, the two main ingredients to be used to prove stability and error estimates are the
ellipticity of the leading diagonal term (here, $(\diff\u,\bv))$, and the {\it inf-sup}
condition. Here the  {\it inf-sup} condition will be  easily provided by the commuting diagram
\eqref{commuting}.  Hence, our main worry will be the treatment of the term
\begin{equation}\label{def:a-mixed}
a(\u,\bv):=(\diff\u,\bv).
\end{equation}
On each element $\E \in \Th$
we define:
\begin{equation}\label{def-ahE}
a^{\E}_h(\bv,\w):=
{{\normalfont{(\diff \br{\bv},\br{\w})_{0,\E}}}}+
S^{\E}(\bv-\br{\bv},\w-\br{\w}),
\end{equation}
where $S^{\E}(\bv,\w)$ is {\it any} symmetric and positive definite bilinear form that {\it scales like} $ a^{\E}(\bv,\w)$ (see \cite{volley}). More precisely, our assumption on $S$ will be: There exist two positive constants {$\alpha_*$ and $\alpha^*$}  (depending on $\diff$  but
independent of $h$) such that
\begin{equation}\label{def:S}
\alpha_* a^{\E}(\bv,\bv)\le S^{\E}(\bv,\bv) \le \alpha^*a^{\E}(\bv,\bv)\quad \forall \bv \in \VMhk.
\end{equation}
For practical purposes it will be convenient to choose the Euclidean scalar product associated to the degrees of freedom in $\VMhk$ multiplied, for instance, by $|\E|\diff({\bf x}_B)$, where ${\bf x}_B=(\xb,\yb)=$  is the barycenter of $\E$.
We notice that, obviously, $\br{{\bf p}}_k={\bf p}_k$ for all ${\bf p}_k\in \Pp_k$. Therefore
\begin{equation}\label{consistency}
a^{\E}_h({\bf p}_k,\w)=\intE \diff  {\bf p}_k\cdot \br{\w},dx \quad \forall\w\in \VMhk, \;\forall{\bf p}_k\in \Pp_k .
\end{equation}

We can now define
\begin{equation}\label{def:ah}
a_h(\bv,\w):=\sum_{\E} a^{\E}_h(\bv,\w).
\end{equation}

\begin{lemma}\label{cont-ellipt-ah}
The bilinear form $a_h(\cdot,\cdot)$ is continuous and elliptic in $(L^2(\Omega))^2$, that is:
\begin{equation}\label{contah}
\begin{aligned}
&\exists M>0 \mbox{ such that } |a_h(\bv,\w)| \le M  \|\bv\|_0 \|\w\|_0 \qquad \forall \bv, \w \in \VMhk,\\
&\exists\, \alpha>0 \mbox{ such that } a_h(\bv,\bv) \ge \alpha \|\bv\|^2_0  \qquad \forall \bv \in \VMhk,
\end{aligned}
\end{equation}
with $M$ and $\alpha$ depending on $\diff$ but independent of $h$.
\end{lemma}
\emph{Proof.} The symmetry of $S^{\E}$ and \eqref{def:S}  imply
easily the continuity of $S^{\E}$:
\begin{equation}\label{contS-0}
S^{\E}(\bv,\w)\le (S^{\E}(\bv,\bv))^{1/2} (S^{\E}(\w,\w))^{1/2}
\le C_{\diff} \|\bv\|_{0,\E}\|\w\|_{0,\E},
\end{equation}
with $C_{\diff}=\alpha^* \diff_{\max}$.
In particular,
\begin{equation}\label{cont:S}
S^{\E}(\bv-\br{\bv},\w-\br{\w})\le C_{\diff} \|\bv-\br{\bv}\|_{0,\E} \|\w-\br{\w}\|_{0,\E}\le C_{\diff}\NErrk(\bv)\NErrk(\w).
\end{equation}
Then, the continuity of $a_h(\cdot,\cdot)$ is an obvious consequence of  the continuity of $a(\cdot,\cdot)$ and of the $L^2-$ projection properties:
\begin{equation*}
 |a_h(\bv,\w)| \le{\diff_{\max}}\|\br{\bv}\|_0\|\br{\w}\|_0+C_{\diff}\|\bv-\br{\bv}\|_0 \|\w-\br{\w}\|_0
 \le M \|\bv\|_0\|\w\|_0.
 \end{equation*}
 Similarly,
 \begin{equation*}
 a_h(\bv,\bv)\ge  \diff_{\min}\Big( \|\br{\bv}\|^2+\alpha_*\|\bv-\br{\bv}\|_0^2\Big)
 \ge \alpha \Big(\|\br{\bv}\|^2_0+\|\bv-\br{\bv}\|^2_0\Big)=\alpha \|\bv\|^2_0.
 \end{equation*}

\hfill\qed

The discrete problem is now:
\begin{equation}\label{vem-approx}
\left\{
\begin{aligned}
&\text{Find~} (\uh,\ph)\in \VMhk\times Q^k_h\text{ such that}\\
&a_h( \uh, \vh)-(\ph,\div \vh)-(\vbeta\cdot\br{\vh}, \ph)=0\quad \forall \vh \in \VMhk\\
&(\div \uh, \qh)+(\reaction \ph,\qh)=(f,\qh)\quad \forall \qh\in Q_h .
\end{aligned}
\right.
\end{equation}
Like we did for the continuous formulation, in order to write \eqref{vem-approx} in a more compact form, we set
$$
\Vh:=\VMhk\times Q^k_h,\quad \UUh:=(\uh,\ph),\quad \VVh:=(\bv,\qh),\quad \FFh:=(0,f),
$$
and
\begin{equation}\label{def:AAh}
\AAh(\UUh,\VVh):=a_h( \uh, \vh)-(\ph,\div \vh)-(\vbeta\cdot \br{\vh}, \ph)
+(\div \uh, \qh)+(\reaction \ph,\qh).
\end{equation}
Then problem \eqref{vem-approx} can be written as
\begin{equation}\label{vem-compact}
\left\{
\begin{aligned}
&\mbox{Find } \UUh\in\Vh \mbox{ such that}\\
&\AAh(\UUh,\VVh)=(\FFh,\VVh)\qquad\forall\,\VVh\in\Vh.
\end{aligned}
\right.
\end{equation}

%-%%%%%%%%%%%%%%%%%%%%--------------------------------------------------------------------
%-%%%%%%%%%%%%%%%%%%%%-
%-%%%%%%%%%%%%%%%%%%%%-
%-%%%%%%%%%%%%%%%%%%%%--------------------------------------------------------------------
%-%%%%%%%%%%%%%%%%%%%%-
%-%%%%%%%%%%%%%%%%%%%%-
\section{Error Estimates}\label{sec:mixed:est}
%-%%%%%%%%%%%%%%%%%%%%--------------------------------------------------------------------
%-%%%%%%%%%%%%%%%%%%%%-
%-%%%%%%%%%%%%%%%%%%%%-
%-%%%%%%%%%%%%%%%%%%%%-
%-%%%%%%%%%%%%%%%%%%%%-
%-%%%%%%%%%%%%%%%%%%%%-

Our final target is to prove the following theorem.
\begin{theorem}\label{risultato}
Under the above assumptions and with the above notation, for $h$ sufficiently small
problem \eqref{vem-approx} has a unique solution $(\uh,\ph)\in \VMhk \times Q^k_h$, and the following error estimates hold:
\begin{equation}\label{stime-finali}
\begin{aligned}
\|\p-\ph\|_0& \le C h^{k+1}\Big(\|\u\|_{k+1}+\|\p\|_{k+1}\Big),\\
 \|\u-\uh\|_0& \le C h^{k+1}\Big(\|\u\|_{k+1}+\|\p\|_{k+1}\Big),\\
 \|\div(\u-\uh)\|_0& \le C h^{k+1}\Big(|f|_{k+1}+\|\p\|_{k+1}\Big),
 \end{aligned}
\end{equation}
with $C$ a constant depending on $\diff, \vbeta,$ and $\reaction$ but independent of $h$.
\end{theorem}
Before proving the theorem, we will introduce some useful lemmata, that deal with properties of the
bilinear forms $\AA$ and $\AAh$.
%
%-%%%%%%%%%%%%%%%%%%%%-
%-%%%%%%%%%%%%%%%%%%%%-
%-%%%%%%%%%%%%%%%%%%%%-
\subsection{Preliminary estimates}
%-%%%%%%%%%%%%%%%%%%%%-
%-%%%%%%%%%%%%%%%%%%%%-
%-%%%%%%%%%%%%%%%%%%%%-
A typical source of difficulties, when proving optimal error estimates, is the fact that the bilinear
form $\AA(\UU,\VV)$ cannot be bounded in terms of the $L^2$ norms of $\UU$ and $\VV$, due to the presence
of the two terms $(\div\u,q)$ and $(p,\div\bv)$ involving the divergence. We will therefore spend some
additional time in order to point out some particular cases in which these terms could be avoided.
In particular, we note that for $\bv\in H(\div;\E)$ and $\q\in L^2(\E)$ we will have
\begin{equation}
\int_{\E}\div\bv\,\q\,\dx =0
\end{equation}
whenever
\begin{itemize}
\item  \hskip0.35truecm $\q\in\Pp_k$, and $\div\bv$ is orthogonal to $\Pp_k$,
\item  $\,\,\,\,\,\div\bv\in\Pp_k$, and $\q$ is orthogonal to $\Pp_k$.
\end{itemize}
Hence, in particular, using \eqref{spazio-VE}, \eqref{spazio-Q}, and \eqref{prop-div} we have:
\begin{equation}
\int_{\E}\div(\w-\Pi^F_h\w)\,\q_h\dx=0\qquad\forall\q_h\in Q^k_h,\quad\forall \w\in (H^1(\E))^2,
\end{equation}
and
\begin{equation}
\int_{\E}\div\vh\,(\r-\Pzk\r)\dx=0\qquad\forall\vh\in \VMhkE,\quad\forall \r\in L^2(\E),
\end{equation}
so that for every $\WW\in\V$ and for every $\VV_h\in\Vh$ we have
\begin{equation}\label{senzadiv}
|\AA(\VV_h,\WW-\WW_I)|+|\AA(\WW-\WW_I,\VV_h)|\le\,C_{\diff,\vbeta,\reaction}\,\|\VV_h\|_0\,\|\WW-\WW_I\|_0.
\end{equation}

%-%%%%%%%%%%%%%%%%%%%%-
%-%%%%%%%%%%%%%%%%%%%%-
%-%%%%%%%%%%%%%%%%%%%%-
\subsection{The consistency error}\label{consi}
%-%%%%%%%%%%%%%%%%%%%%-
%-%%%%%%%%%%%%%%%%%%%%-
%-%%%%%%%%%%%%%%%%%%%%-
%%
Further attention should also be given to the difference $(\AA_h-\AA)(\WW,\VV)$. We will perform the analysis on a single element, without indicating every time that the norms are considered
in $L^2(\E)$. Using \eqref{def:AA0} and \eqref{def:AAh} we have easily
\begin{equation}\label{termine-III}
\begin{aligned}
(\AA_h-\AA)(\WW,\VV)&=(\diff \br \w,\br\bv)-(\diff \w,\bv)~(=:\text{T}_1(\WW,\VV))\\
&+S(\w-\br{\w},\bv-\br{\bv})~(=:\text{T}_2(\WW,\VV))\\
&+(\br{\bv}-\bv,\vbeta\r)~(=:\text{T}_3(\WW,\VV)),
\end{aligned}
\end{equation}
where as before $\VV=(\bv,q)$ and $\WW=(\w,r)$ are in $\Vh$.
 We point out that all the terms $\text{T}_1$, $\text{T}_2$ and $\text{T}_3$ do not involve derivatives, so
 that we will not have continuity problems. For the term $\text{T}_1$, using repeatedly the properties of the $L^2-$ projection we have:
 \begin{equation}\label{consi1}
\begin{aligned}
\text{T}_1(\WW,\VV)&=(\diff \br \w,\br\bv)-(\diff \w,\bv)=(\diff \w,\br\bv-\bv)-( \w-\br\w,\diff\br \bv)\\
&=(\diff \w-\overline{\diff\w},\br\bv-\bv)-( \w-\br\w,\diff\br \bv-\overline{\diff \bv})\\
&=(\diff \w-\overline{\diff\w},\br\bv-\bv)-(\w-\br \w,\diff\br\bv-\overline{\diff \bv}+\diff\bv-\diff\bv)\\
&=(\diff \w-\overline{\diff\w},\br\bv-\bv)-(\w-\br \w,\diff\bv-\overline{\diff \bv})-(\w-\br\w,\diff(\br\bv-\bv))\\
&\le\Big(C_{\diff}\|\w-\br{\w}\|_0 +\|\diff\w-\br{\diff\w}\|_0\Big)\|\bv\|_0 \\
&\le \Big(C_{\diff}\NErrk(\WW)+\NErrk(\diff\WW)\Big)\|\VV\|_0 .
\end{aligned}
\end{equation}
Needless to say, in view of the symmetry of the term, we also have
\begin{equation}\label{consi2}
\text{T}_1(\WW,\VV)=\text{T}_1(\VV,\WW)\le \,\Big(C_{\diff}\NErrk(\VV)+\NErrk(\diff\VV)\Big)\|\WW\|_0 .
\end{equation}
The terms $\text{T}_2$ and $\text{T}_3$ in \eqref{termine-III} are easily bounded.
Directly from \eqref{cont:S} we have
\begin{equation}\label{consi5}
\text{T}_2(\WW,\VV)\le C_{\diff} \NErrk(\WW)\,\NErrk(\VV),
\end{equation}
and for $\text{T}_3$
\begin{equation}\label{consi4}
\text{T}_3(\WW,\VV)=(\br{\bv}-\bv,\vbeta\r-\br{\vbeta\r})\le
\,\NErrk(\VV)\, \NErrk(\vbeta\WW).
\end{equation}

\medskip \noindent  The above analysis can now be summarized in the following two estimates, that will both be used in our final proof.
\begin{itemize}
\item Using \eqref{consi2}, \eqref{consi4} with \eqref{consi6}, and \eqref{consi5} we have
\begin{equation}\label{consi7}
(\AA_h-\AA)(\WW,\VV)\le \, C_{\diff,\vbeta}\,\Big((\NErrk(\VV)+\NErrk(\diff\VV)\Big)\|\WW\|_0 .
\end{equation}
\item Using instead \eqref{consi1}, \eqref{consi4} and \eqref{consi5} we have
\begin{equation}\label{consi8}
\!(\AA_h-\AA)(\WW,\VV)\le \,C_{\diff}\,\Big(\NErrk(\WW)
+\NErrk(\diff\WW)+\NErrk(\vbeta\WW)\Big)\|\VV\|_0 .
\end{equation}
\end{itemize}

%-%%%%%%%%%%%%%%%%%%%%-
%-%%%%%%%%%%%%%%%%%%%%-
%-%%%%%%%%%%%%%%%%%%%%-
\subsection{The dual problem}
%-%%%%%%%%%%%%%%%%%%%%-
%-%%%%%%%%%%%%%%%%%%%%-
%-%%%%%%%%%%%%%%%%%%%%

Our proof will use a duality argument. Therefore we spend some time analyzing the dual problem.

\begin{lemma}\label{dual} Let $\fd\in L^2(\Omega)$, $\bg\in H(\div;\Omega)$, and set ${\bf G}:=(\bg ,\fd)$. Let
$\ZZ:=(\bzeta,z) \in \V$ be the solution of
\begin{equation}\label{eqzeta-2}
\AA(\WW,\ZZ)=({\bf G},\WW)\qquad \forall\, \WW = (\w,\r)\in\V.
\end{equation}
Then $\ZZ$ is the  solution of
\begin{equation}\label{splitted-d}
\bzeta=\diffp(\nabla \z+\bg)~~\text{ and }~~ -\div\bzeta-\vbeta\cdot\bzeta + \reaction \,\z=\fd~~\text{in }
 \Omega,\quad \z=0 ~~\text{ on }\Gamma
\end{equation}
that is, (see \eqref{def:adjoint}),
\begin{equation}\label{eq:z}
\Amfs \z
=\fd+\bb\cdot\bg+ \div(\diffp\bg),
\end{equation}
so that, in particular
\begin{equation}\label{ellreg}
||\z||_2+||\bzeta||_1\le C^*(||\fd||_0+||\diffp\bg||_{H(\div)}).
\end{equation}
\end{lemma}
\emph{Proof.} Recalling \eqref{def:AA0}, and substituting $\WW$
for $\UU$ and $\ZZ$ for $\VV$ we get
\begin{equation}\label{def:AA-3}
\AA(\WW,\ZZ)=(\diff\w, \bzeta)-(\r,\div \bzeta)-(\vbeta\cdot \bzeta, \r)
+(\div \w, \z)+(\reaction \r,\z) .
\end{equation}
Separating the equations in $\w$ and in $\r$ in \eqref{eqzeta-2} it is not difficult to see that $(\bzeta, z)$ solves
\begin{equation}\label{mixed-conts}
\left\{
\begin{aligned}
&(\diff\w, \bzeta)+(\div \w, \z)=(\bg,\w)\quad \forall \w \in H(\div,\Omega)\\
&-(\r,\div \bzeta)-(\vbeta\cdot\bzeta, \r)+(\reaction \r,\z)=(\fd,\r)\quad \forall r\in L^2(\Omega)
\end{aligned}
\right.
\end{equation}
giving, respectively, $$\bzeta=\diffp\nabla\z+\diffp\bg\quad\mbox{plus } \z\in H^1_0(\Omega),$$ and $$-\div\bzeta-\vbeta\cdot\bzeta+\reaction\z=\fd.$$
Putting them together we have
$$-\div(\diffp\nabla\z)-\div(\diffp\bg)-\bb\cdot\nabla\z-\bb\cdot\bg+\reaction\z=\fd ,$$
and \eqref{eq:z} follows.

\qed

%-%%%%%%%%%%%%%%%%%%%%-
%-%%%%%%%%%%%%%%%%%%%%-
%-%%%%%%%%%%%%%%%%%%%%-
\subsection{Proof of Theorem \ref{risultato}}
%-%%%%%%%%%%%%%%%%%%%%-
%-%%%%%%%%%%%%%%%%%%%%-
%-%%%%%%%%%%%%%%%%%%%%-

We are now ready for the proof of Theorem \ref{risultato}.

\emph{Proof.} To prove Theorem \ref{risultato} we shall follow the
arguments of Douglas-Roberts \cite{DR-85}. We first assume that
\eqref{vem-approx} has a solution, at least for $h$ sufficiently
small. That it does, it will be clear from the convergence
analysis. Let therefore $\UU_h= (\uh,\ph)$ be a solution of
\eqref{vem-approx}. Let us form the error equation:
\begin{equation}\label{error-equations-compact}
\AA(\UU,\VV_h)-\AA_h(\UU_h,\VV_h)=0\;\quad \forall\,\VV_h\equiv (\bv_h,q_h)\in\Vh.
\end{equation}
We use duality arguments. Let $\bpsi=(\bchi,\psi)$ be the solution of the adjoint problem
\begin{equation}\label{adjoint-compact}
\AA(\VV,\bpsi)=\Big(\diff(\UU_I-\UU_h),\VV\Big)=\Big((\diff(\u_I-\u_h),p_I-p_h),\VV\Big) \quad \forall \VV \in \V.
\end{equation}
According to Lemma \ref{dual},  $\psi\in H^1_0(\Omega)\cap H^2(\Omega)$ is the solution of the adjoint problem
\begin{equation}\label{eq:psi}
\Amfs \psi\equiv\div(-\diffp\nabla \psi)- \bb\cdot\nabla \psi + \reaction\,\psi=\pI-\ph+\vbeta\cdot (\uI-\uh)+ \div(\uI-\uh),
\end{equation}
and by the elliptic regularity \eqref{ellreg} with ${\bf G}\equiv (\bg,\fd):=(\diff(\u_I-\u_h),p_I-p_h)$
we get
\begin{equation}\label{stimapsi1}
\|\psi\|_2+\|\bchi\|_1\le C^*\, (\|\pI-\ph\|_0 +\|\uI-\uh\|_{H(\div)}) .
\end{equation}
Our first step will then be the estimate of $\|\div(\uI-\uh)\|_0$.
Looking at the discrete and continuous  equations we have
\begin{equation}\label{spaccapalle}
\div\uh=\Pzk(f-\reaction \ph ) \quad \mbox{and }\quad \div\u=f-\reaction \p,
\end{equation}
and from \eqref{prop-div}  $\div\uI=\Pzk\div\u=\Pzk(f-\reaction \p)$. Hence,
\begin{equation}
 \div(\uI-\uh)=\Pzk(\reaction (\ph-\p)),
\end{equation}
so that, clearly,
\begin{equation}
\|\div(\uI-\uh)\|_0\le C_{\reaction}\|\p-\ph\|_0.
\end{equation}
Therefore, \eqref{stimapsi1} reduces to
\begin{equation}\label{stimapsi2}
\begin{aligned}
\|\psi\|_2+\|\bchi\|_1&\le C\,( \|\pI-\ph\|_0 + \|\uI-\uh\|_0 + \|p-p_I\|_0)\\
&\le C \Big(\|\UU_I-\UU_h\|_0 + \NErrk(\UU)\Big),
\end{aligned}
\end{equation}
and using (for instance) \eqref{inH1} with $\wiz=1$ the estimate \eqref{stimapsi2}  implies that
\begin{equation}\label{stimapsi3}
\NErrk(\bpsi)\le C\,h\, \Big(\|\psi\|_1 + \|\bchi\|_1\Big)\le C\,h\, \Big(\|\UU_I-\UU_h\|_0 + \NErrk(\UU)\Big),
\end{equation}
as well as
\begin{equation}\label{stimapsi31}
\|\bpsi_I\|_0\le  \|\bpsi-\bpsi_I\|_0+\|\bpsi\|_0 \le C \Big( \|\UU_I-\UU_h\|_0 + \NErrk(\UU)\Big).
\end{equation}
Moreover, taking $\VV=\UU_I-\UU_h$ in \eqref{adjoint-compact}, it is immediate to see that
\begin{equation}
{\diff}_ {\min}\|\UU_I-\UU_h\|^2\le \int_{\Omega}( \diff|\u_I-\u_h|^2+|p_I-p_h|^2) \dx
=\AA(\UU_I-\UU_h,\bpsi).
\end{equation}
Hence,
\begin{equation}\label{stima1}
\begin{aligned}
{\diff}_{\min} &\|\UU_I-\UU_h\|^2\le \,\AA(\UU_I-\UU_h,\bpsi)\;(\pm \bpsi_I)\\
&=\AA(\UU_I-\UU_h,\bpsi-\bpsi_I)+\AA(\UU_I-\UU_h,\bpsi_I)\;(\pm \UU)\\
&=\text{I}+\AA(\UU_I-\UU,\bpsi_I)+\AA(\UU-\UU_h,\bpsi_I)\;(\mbox{just linearity})\\
&=\text{I}+\text{II}+\AA(\UU,\bpsi_I)-\AA(\UU_h,\bpsi_I) \,\mbox{ (use \eqref{error-equations-compact})}\\
&=\text{I}+\text{II}+(\AA_h-\AA)(\UU_h,\bpsi_I).
\end{aligned}
\end{equation}
The first two terms are easily bounded using \eqref{senzadiv}, \eqref{stimapsi3}-\eqref{stimapsi31}, and \eqref{stime-variE}:
\begin{equation}\label{termine-I}
\begin{aligned}
\text{I}\equiv\AA(\UU_I-\UU_h,\bpsi-\bpsi_I)&\le
C\,\|\UU_I-\UU_h\|_0\, h\Big(\|\UU_I-\UU_h\|_0+\NErrk(\UU)\Big)\\
&\le C \Big(h \|\UU_I-\UU_h\|^2_0 + h^{k+2} \|\UU_I-\UU_h\|_0\Big),
\end{aligned}
\end{equation}
\begin{equation}\label{termine-II}
\begin{aligned}
\text{II}\equiv\AA(\UU_I-\UU,\bpsi_I)&\le C\, \NErrk(\UU)\,\Big(\|\UU_I-\UU_h\|_0+\NErrk(\UU)\Big)\\
&\le C\, \Big(\|\UU_I-\UU_h\|_0 h^{k+1}+h^{2k+2}\Big),
\end{aligned}
\end{equation}
and we are left with the third term. For it, we are going to use the arguments of Subsection \ref{consi}.
We start by observing that
\begin{equation}\label{termine-III-0}
(\AA_h-\AA)(\UU_h,\bpsi_I)=(\AA_h-\AA)(\UU_h-\UU_I,\bpsi_I)+(\AA_h-\AA)(\UU_I,\bpsi_I).
\end{equation}
The first term in \eqref{termine-III-0} can be easily bounded, using \eqref{consi7},  \eqref{errWI}, \eqref{selava},  \eqref{stimapsi3}, and \eqref{stime-variE}:
\begin{equation}\label{termine-III-1}
\begin{aligned}
(\AA_h-\AA)(\UU_h-\UU_I,\bpsi_I)&\le
 C_{\diff,\vbeta}\Big(\NErrk(\bpsi_I)+\NErrk(\diff\bpsi_I)\Big)\|\UU_h-\UU_I\|_0\\
 &\le C\,h\,\Big(\|\UU_h-\UU_I\|_0+\NErrk(\UU)\Big)\,\|\UU_h-\UU_I\|_0\\
 &\le C\,\Big(h\,\|\UU_h-\UU_I\|^2_0+h^{k+2} \|\UU_h-\UU_I\|_0\Big),
\end{aligned}
\end{equation}
while, using \eqref{consi8}, \eqref{errWI}, \eqref{selava}, and \eqref{stimapsi31}, the second term in \eqref{termine-III-0} can be bounded
by
\begin{equation}\label{termine-III-2}
\begin{aligned}
(\AA_h-\AA)&(\UU_I,\bpsi_I)\le
\,C_{\diff}\,\Big(\NErrk(\UU_I)+\NErrk(\diff\UU_I)+\NErrk(\vbeta\UU_I)\Big)\|\bpsi_I\|_0\\
&\le C\,\Big(\NErrk(\UU)
+\NErrk(\diff\UU))
+\NErrk(\vbeta\UU)\Big)
\Big(\|\UU_h-\UU_I\|_0+\NErrk(\UU)\Big)\\
&\le C\,\Big(h^{k+1} \|\UU_h-\UU_I\|_0+ h^{2k+2}
\Big).
\end{aligned}
\end{equation}
Inserting \eqref{termine-I}, \eqref{termine-II}, and \eqref{termine-III-1}-\eqref{termine-III-2} into
\eqref{stima1}
we have then
\begin{equation}\label{pen}
 \diff_{\min} \|\UU_h-\UU_I\|^2
\le C\,\Big( h\|\UU_h-\UU_I\|^2+\|\UU_h-\UU_I\|\,h^{k+1}+h^{2k+2}\Big).
\end{equation}
For $h$ small enough (say: $C h\le(1/2)\diff_{min}$ in \eqref{pen}) we can hide the first term in the r.h.s. of \eqref{pen} in the left-hand side, and have
\begin{equation}
 \|\UU_h-\UU_I\|_0^2\le\,C\, \Big(h^{k+1} \|\UU_h-\UU_I\|_0\,+h^{2k+2} \Big),
\end{equation}
and the first two estimates in \eqref{stime-finali} follow
completing the square. The estimate on the divergence follows
directly from \eqref{spaccapalle},  and standard error estimates.

Finally, since \eqref{vem-approx} is finite dimensional, {in order
to prove the existence of the solution} we only have to prove
uniqueness, that is, we have to prove that for $f=0$ problem
\eqref{vem-approx} has only the solution $\ph=0, \uh=0$. Since we
assumed that the continuous problem \eqref{Pb-cont} has a unique
solution, it follows that for $f=0$ we have $\p=0, \u=0$. The
above analysis showed that, for $h$ small enough, any solution
$(\uh,\ph)$ of \eqref{vem-approx} must satisfy
\eqref{stime-finali} which, in our case, imply $\uh=0,\ph=0$, and
the proof is concluded.

\qed

\begin{remark}\label{2-3Dim}Looking at the construction of the method, and to the analysis of its convergence properties, it is not difficult to see that
%, in the present case where only VEM approximations of $H(\div;\Omega)$ and polynomial approximations of $L^2(\Omega)$ are used,
the passage from the two-dimensional case to the three-dimensional one can be done, using \cite{super-misti}, without any difficulty. However, the notation for dealing with both cases at the same time would be more cumbersome, and a presentation with two separate treatments would be very boring and essentially useless.
\end{remark}

%-%%%%%%%%%%%%%%%%%%%%-
%-%%%%%%%%%%%%%%%%%%%%-
%-%%%%%%%%%%%%%%%%%%%%-
\section{Superconvergence results}\label{superconvergence}
%-%%%%%%%%%%%%%%%%%%%%-
%-%%%%%%%%%%%%%%%%%%%%-
%-%%%%%%%%%%%%%%%%%%%%-

\begin{theorem}
Let $\ph$ be the solution of \eqref{vem-approx}, and let $\pI\in Q^k_h$ be the interpolant of $p$. Then, for $h$ sufficiently small,
\begin{equation}\label{superconv}
\|\pI-\ph\|_0 \le C\, h^{k+2} \Big(\|\u\|_{k+1}+\|\p\|_{k+1}+|f|_{k+1}\Big),
\end{equation}
where $C$ is a constant depending on $\diff, ~\vbeta,$ and $\reaction$ but independent of $h$.
\end{theorem}
\emph{Proof.} We proceed again via duality argument. Let $\psi \in
H^1_0(\Omega)\cap H^2(\Omega)$ be the solution of the adjoint
problem
\begin{equation}\label{def:adjoint}
\div(-\diffp(\xx)\nabla\psi){- \bb(\xx)\cdot\nabla\psi} {+ \reaction(\xx)\,\psi}= \pI-\ph,\qquad \bchi=\diffp\nabla\psi,
\end{equation}
whose mixed formulation is{: Find $(\bchi,\psi)$ in $H(\div,\Omega)\times  L^2(\Omega)$  such that}
\begin{equation}\label{mixed-adjoint}
\left\{
\begin{aligned}
&(\diff \bchi, \bv)+(\psi,\div \bv)= 0 {\quad \forall \bv\in H(\div,\Omega)}\\
&-(\div \bchi, q){-(\vbeta\cdot \bchi, q)}{+(\reaction \psi,q)}=(\pI-\ph,q)
{\quad \forall \q\in L^2(\Omega)}.
\end{aligned}
\right.
\end{equation}
The error equations \eqref{error-equations-compact}, using \eqref{medie} and \eqref{prop-div}, become
\begin{equation}\label{error-eqn}
\left\{
\begin{aligned}
&a(\u,\vh)-a_h( \uh, \vh)-(\pI-\ph,\div \vh){-(\vbeta\cdot \vh, \p)+(\vbeta\cdot\br{\vh}, \ph)}=0\quad \forall \vh \in \VMhk ,\\
&(\div( \u-\uh), \qh){+(\reaction (\p- \ph),\qh)}=0\quad \forall \qh\in Q_h .
\end{aligned}
\right.
\end{equation}
Taking now $q=\pI-\ph$ in \eqref{mixed-adjoint} gives
\begin{equation}\label{def:norma0}
\|\pI-\ph\|^2_0= -(\div \bchi,\pI-\ph){-(\vbeta\cdot \bchi, \pI-\ph)}{+(\reaction \psi,\pI-\ph)} .
\end{equation}
For the first term, using again \eqref{medie} and \eqref{prop-div}, we have
%{ho aumentato il rientro equazione seguente}
\begin{equation}\label{stimaA}
\begin{aligned}
(\div \bchi,\pI&-\ph)= (\div \bchi_I,\pI-\ph)\qquad (use \eqref{error-eqn} \mbox{ with }\vh=\bchi_I)\\
&=a(\u,\bchi_I)-a_h(\uh,\bchi_I){-(\vbeta\cdot \bchi_I, \p)+(\vbeta\cdot\br{\bchi_I}, \ph)}\qquad(\pm \uh)\\
&=a(\u-\uh,\bchi_I)+a(\uh,\bchi_I)-a_h(\uh,\bchi_I)
{-(\vbeta\cdot \bchi_I, \p)+(\vbeta\cdot\br{\bchi_I}, \ph)}\qquad(\pm \bchi)\\
&=a(\u-\uh,\bchi)+a(\u-\uh,\bchi_I-\bchi)+a(\uh,\bchi_I)-a_h(\uh,\bchi_I)\\
&{-(\vbeta\cdot \bchi_I, \p)+(\vbeta\cdot\br{\bchi_I}, \ph)} .
\end{aligned}
\end{equation}
In turn, the first term in \eqref{stimaA} becomes
\begin{equation}\label{first-super}
\begin{aligned}
a(\u-\uh,\bchi)&=(\u-\uh, \nabla \psi)=-(\div(\u-\uh), \psi)\qquad(\pm \psi_I)\\
&=-(\div(\u-\uh), \psi-\psi_I)-(\div(\u-\uh), \psi_I)\qquad(\mbox{use }\eqref{error-eqn})\\
&=-(\div(\u-\uh), \psi-\psi_I) {+(\reaction \psi_I,\p-\ph)} .
%&\le C h \|\div(\u-\uh)\|_0 \|\pI-\ph\|_0
\end{aligned}
\end{equation}
Replacing \eqref{first-super} in \eqref{stimaA}, and using the result for the first term of \eqref{def:norma0}, we have then
\begin{equation}\label{def:norma0bis}
\begin{aligned}
\|\pI-\ph\|^2_0&=-\big[-(\div(\u-\uh), \psi-\psi_I) {+(\reaction \psi_I,\p-\ph)} +a(\u-\uh,\bchi_I-\bchi)\\
&+a(\uh,\bchi_I)-a_h(\uh,\bchi_I)
{-(\vbeta\cdot \bchi_I, \p)+(\vbeta\cdot\br{\bchi_I}, \ph)}\big]\\
&{-(\vbeta\cdot \bchi, \pI-\ph)}{+(\reaction \psi,\pI-\ph)}.
\end{aligned}
\end{equation}
The first two terms are easily bounded:
\begin{equation}\label{stima-diff}
\begin{aligned}
|a(\u-\uh,\bchi_I-\bchi)|&\le C h \|\u-\uh\|_0 \|\pI-\ph\|_0,\\
|(\div(\u-\uh), \psi-\psi_I)|&\le C h^2 \|\div(\u-\uh)\|_0 \|\pI-\ph\|_0,
\end{aligned}
\end{equation}
while using \eqref{consi2} and \eqref{consi5} we get
\begin{equation}
|a(\uh,\bchi_I)-a_h(\uh,\bchi_I)|\le C_{\nu} h^{k+1} \|\u\|_{k+1,\Omega} h\, \|\pI-\ph\|_0.
\end{equation}
For the terms involving reaction, adding and subtracting ${(\reaction \psi_I, \pI-\ph)}$ and using the properties of the projection we obtain
\begin{equation}\label{stima-reaction}
\begin{aligned}
{(\reaction \psi,\pI-\ph)} {-(\reaction \psi_I,\p-\ph)}&=
%\qquad(\pm \psi_I)\\
{(\reaction (\psi-\psi_I),\pI-\ph)+(\reaction \psi_I,\p_I-\ph)-(\reaction \psi_I,\p-\ph)}\\
&={(\reaction (\psi-\psi_I),\pI-\ph)+(\reaction \psi_I,\p_I-\p)}\\
&={(\reaction (\psi-\psi_I),\pI-\ph)+(\reaction \psi_I-\br{\reaction \psi_I},\p_I-\p)}\\
&\le {C_{\reaction} h^2 \Big(\|\pI-\ph\|^2_0 + \|\p-\p_I \| \|\pI-\ph\|_0\Big)}.
\end{aligned}
\end{equation}
For $h$ small enough the first term in the right-hand side of \eqref{stima-reaction} can be hidden in the left-hand side of \eqref{def:norma0bis} and the other one is more than enough.

\noindent Finally, the terms involving advection can be treated as:
\begin{equation}\label{stima-advection}
\begin{aligned}
&-(\vbeta\cdot \bchi, \pI-\ph)+(\vbeta\cdot \bchi_I, \p)-(\vbeta\cdot\br{\bchi_I}, \ph) \qquad(\pm \bchi_I)\\
&=-(\vbeta\cdot (\bchi-\bchi_I), \pI-\ph)-(\vbeta\cdot \bchi_I, \pI-\ph)+(\vbeta\cdot \bchi_I, \p)-(\vbeta\cdot\br{\bchi_I}, \ph)\\
&=-(\vbeta\cdot (\bchi-\bchi_I), \pI-\ph)+(\vbeta\cdot \bchi_I, \p-\pI)+(\bchi_I-\br{\bchi_I}, \vbeta\ph)\\
&=-(\vbeta\cdot (\bchi-\bchi_I), \pI-\ph)+(\vbeta\cdot \bchi_I-\br{\vbeta\cdot \bchi_I}, \p-\pI)+(\bchi_I-\br{\bchi_I}, \vbeta\ph-\br{\vbeta\ph})\\
&\le C_{\vbeta} h \|\pI-\ph\|_0 \Big( \|\p-\p_I \|_0 +\|\p-\ph\|_0+ h^{k+1}\|\p\|_{k+1}\Big).
\end{aligned}
\end{equation}
Inserting \eqref{stima-diff}--\eqref{stima-advection} in \eqref{def:norma0bis} and using \eqref{stime-finali} and standard interpolation estimates we obtain \eqref{superconv}
%The first term, for $h$ small enough can be hidden in the left-hand side of \eqref{def:norma0bis}, and the other two are good. For example:
%\begin{equation}
%{\color{red}
%(\br{\bchi_I}- \bchi_I, \vbeta\ph-\br{\vbeta\ph})=(\br{\bchi_I}- \bchi_I, \vbeta\ph-\vbeta \p+\vbeta \p-\br{\vbeta \p}+\br{\vbeta \p}-\br{\vbeta\ph})}
%\end{equation}

\qed

\section{Numerical Experiments}
\label{Num-Exp}

In this Section we will present some numerical experiments to validate the convergence
results proven in the previous sections. We will test our method on the same problem
and with the same meshes of \cite{variable-primal}, where we studied the Virtual
Element Method for problem \eqref{Pb-cont} in the primal form.

Before presenting the numerical results we make a comment on the stabilization bilinear form in \eqref{def-ahE}.
%-------------------------------------------------------------------------------------------------------------------------
For each element $\E \in \Th$ we denote by $\chi_i$, for $i=1,2,..,N_\E$, the operator $\VMhkE \rightarrow {\mathbb R}$ that to each $\vh\in\VMhkE$ associates the $i$-th local degree of freedom \eqref{dof1}-\eqref{dof2}-\eqref{dof3}, ordered as follows: first the boundary d.o.f.  \eqref{dof1}, for $i=1,2,...,N_\E^\partial$, and then the internal ones \eqref{dof2}-\eqref{dof3}, for $i=N_\E^\partial+1,...,N_\E$.
%We moreover order such operators in the following way: for $i=1,2,...,N_\E^\partial$ the $\chi_i$ are associated to the boundary degrees of freedom \eqref{dof1} while for $i=N_\E^\partial+1,...,N_\E$ the $\chi_i$ are associated to the internal ones \eqref{dof2}-\eqref{dof3}.
%For each element $\E \in \Th$ we denote by $\chi_i$, for $i=1,2,..,N_\E$, the operator $\VMhkE \rightarrow {\mathbb R}$ that to each $v_h\in\VMhkE$ associates the i-th local degree of freedom \eqref{dof1}-\eqref{dof2}-\eqref{dof3}, and we order first the
%We moreover order such operators in the following way: for $i=1,2,...,N_\E^\partial$ the $\chi_i$ are associated to the boundary degrees of freedom \eqref{dof1} while for $i=N_\E^\partial+1,...,N_\E$ the $\chi_i$ are associated to the internal ones \eqref{dof2}-\eqref{dof3}.
%%
We assume that all the degrees of freedom are scaled in such a way that the associated dual basis $\{ \phi_i \}_{i=1}^{N_E}$ scales uniformly in the mesh size
\begin{equation}
|| \phi_i ||_{L^\infty(E)} \simeq 1 \quad \forall i=1,2,...,N_E .
\end{equation}
With this notation, the most natural VEM stabilization $S^\E(\cdot,\cdot)$ in \eqref{def-ahE} is given by (see \cite{volley})
\begin{equation}\label{L-new-1}
S^{\E}(\bv-\PPzk \bv , \w-\PPzk\w) := |\E| \sum_{i=1}^{N_{\E}} \chi_i \big( \bv-\PPzk\bv \big) \: \chi_i\big(\bw-\PPzk\w\big)
\end{equation}
for all $\bv,\w\in\VMhkE$.
We now observe that, by definition of the $L^2$ projection operator $\PPzk$, and since both spaces $\calG_{k-1}(\E),\, \calG^\perp_{k}(\E)$ appearing in \eqref{dof2}-\eqref{dof3} are included in $(\Pp_k(\E))^2$, it is immediate to check that
$$
\chi_i \big( \bv-\PPzk\bv \big) = 0 \quad \forall \bv\in\VMhkE, \ i=N_\E^\partial+1,...,N_\E .
$$
Therefore the contribution of the internal degrees of freedom in \eqref{L-new-1} vanishes, and we can equivalently use the shorter version
$$
S^{\E}(\bv-\PPzk \bv , \w-\PPzk\w) := |\E| \sum_{i=1}^{N_{\E}^\partial} \chi_i \big( \bv-\PPzk\bv \big) \: \chi_i\big(\bw-\PPzk\w\big).
$$
In other words, the internal degrees of freedom do not need to be included in the stabilization procedure.

\newcommand{\pex}{p}
\newcommand{\uex}{\boldsymbol{u}}

\newcommand{\xmax}{x_\textup{max}}
\newcommand{\ymax}{y_\textup{max}}

\subsection{Exact Solution}
We will consider problem \eqref{Pb-cont} on the unit square with
\begin{equation}
\diffp(x,y) = \begin{pmatrix}y^2+1&-xy\\-xy&x^2+1\end{pmatrix},\quad
\bb=(x,y),\quad
\reaction = x^2+y^3,
\end{equation}
and with right hand side and Dirichlet boundary conditions defined
in such a way that the exact solution is
\begin{equation}
\pex(x,y) = x^2 y + \sin(2\pi x) \sin(2\pi y)+2.
\end{equation}
The corresponding flux is given by
\begin{equation}
\uex = -\diffp\nabla\pex + \boldsymbol{b}\,\pex.
\end{equation}
We will show, in a loglog scale, the convergence curves of the error in $L^2$
between $(\pex,\uex)$ and the solution $(\ph,\uh)$
given by the mixed Virtual Element Method \eqref{vem-approx}.
As the VEM flux $\uh$ is not explicitly known inside the elements, we compare $\uex$ with the $L^2-$projection of $\uh$ onto $(\Pp_k)^2$, that is, with $\PPzk\uh$.
%We will also show the behaviour of $|\pex -\Pi^0_k\,\p_h|$ at the
%maximum point of $\pex$ which is approximately at $(\xmax,\ymax)=(0.781,0.766)$.
%which is internal to an element in all the meshes.
%

\subsection{Meshes}

For the convergence test we consider four sequences of meshes.% described below.

\begin{figure}
\hfill
 \begin{minipage}[b]{0.45\textwidth}
  \begin{center}
  \includegraphics[width=\textwidth]{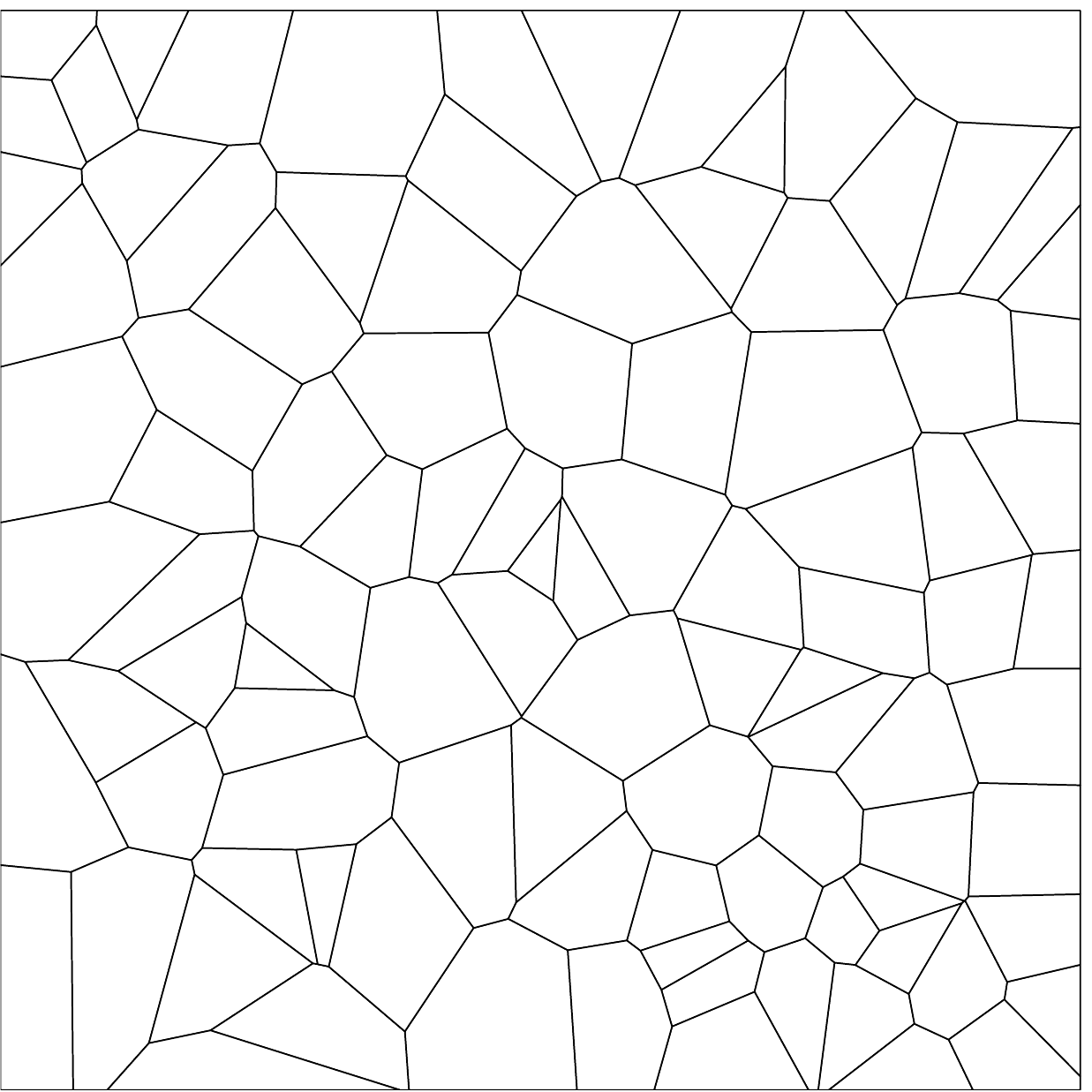}
  \end{center}
  \caption{\texttt{Lloyd-0} mesh}
 \label{fig:Lloyd-0}
 \end{minipage}
\hfill
 \begin{minipage}[b]{0.45\textwidth}
  \begin{center}
  \includegraphics[width=\textwidth]{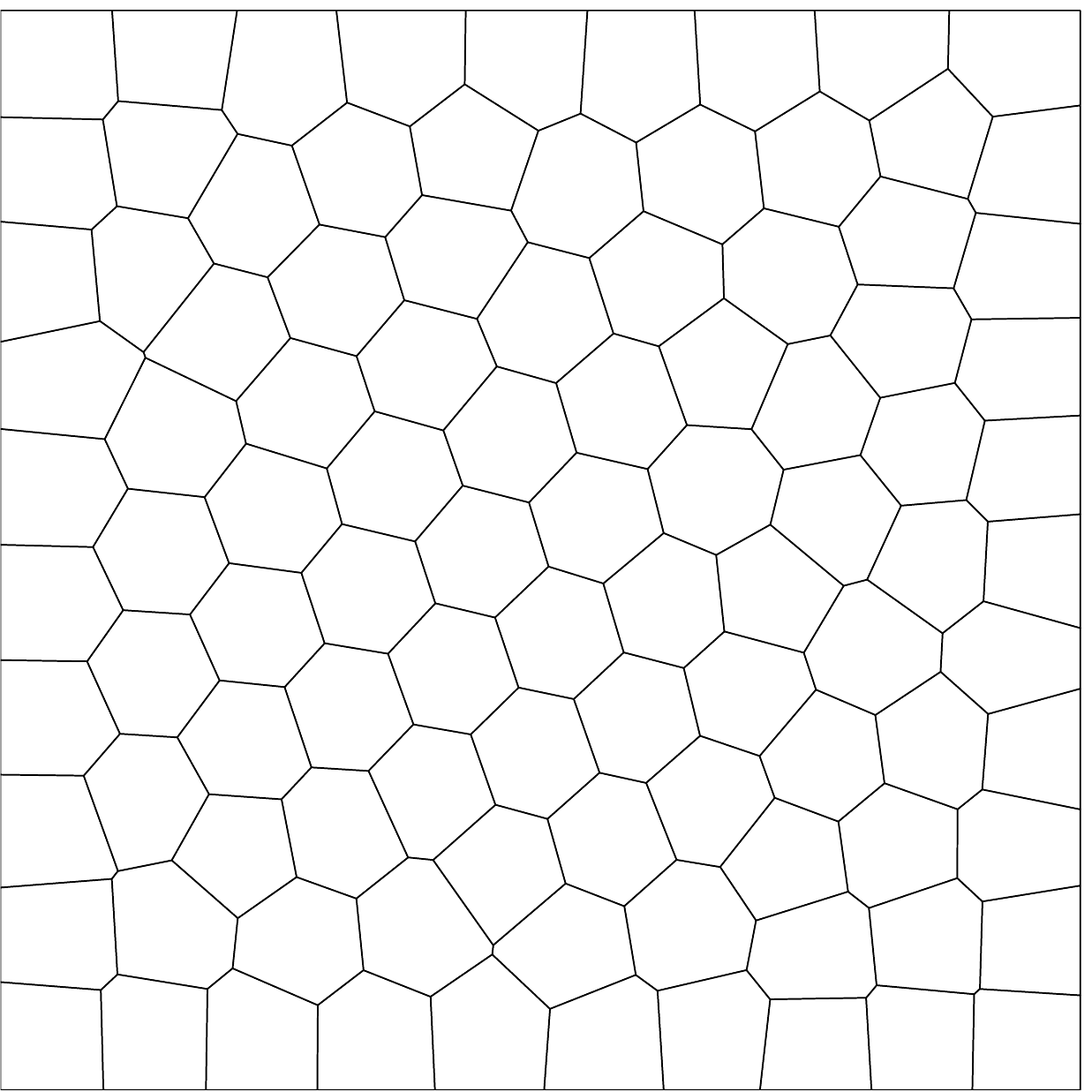}
  \end{center}
  \caption{\texttt{Lloyd-100} mesh}
 \label{fig:Lloyd-100}
 \end{minipage}
\hfill
\end{figure}

The first sequence of meshes (labelled \texttt{Lloyd-0}) is a
random Voronoi polygonal tessellation of the unit square
in 25, 100, 400 and 1600 polygons.
The second sequence (labelled \texttt{Lloyd-100})
is obtained starting from the previous one and
performing 100 Lloyd iterations  leading to a Centroidal Voronoi Tessellation (CVT)
(see e.g.  \cite{Du:Faber99}).
%
%The first two sequences of meshes are obtained from the
%PolyMesher software \cite{TPPM12}.
%PolyMesher is a polygonal mesh generator that starts from a
%random Voronoi mesh and then regularizes it by means of Lloyd iterations,
%leading to a Centroidal Voronoi Tessellation (CVT).
%%
%The first sequence of meshes (labelled \texttt{Lloyd-0}) is just the first try of
%PolyMesher, i.e., a random Voronoi tessellation of the square, while the second family
%(labelled \texttt{Lloyd-100}) is obtained after 100 Lloyd iterations
%and is much more regular. Both sequences consist of meshes with
%25, 100, 400 and 1600 polygons.
The 100-polygon mesh of each family is shown in Fig.~\ref{fig:Lloyd-0}
(\texttt{Lloyd-0}) and in Fig.~\ref{fig:Lloyd-100} (\texttt{Lloyd-100})
respectively.

\begin{figure}
\hfill
 \begin{minipage}[b]{0.45\textwidth}
  \begin{center}
  \includegraphics[width=\textwidth]{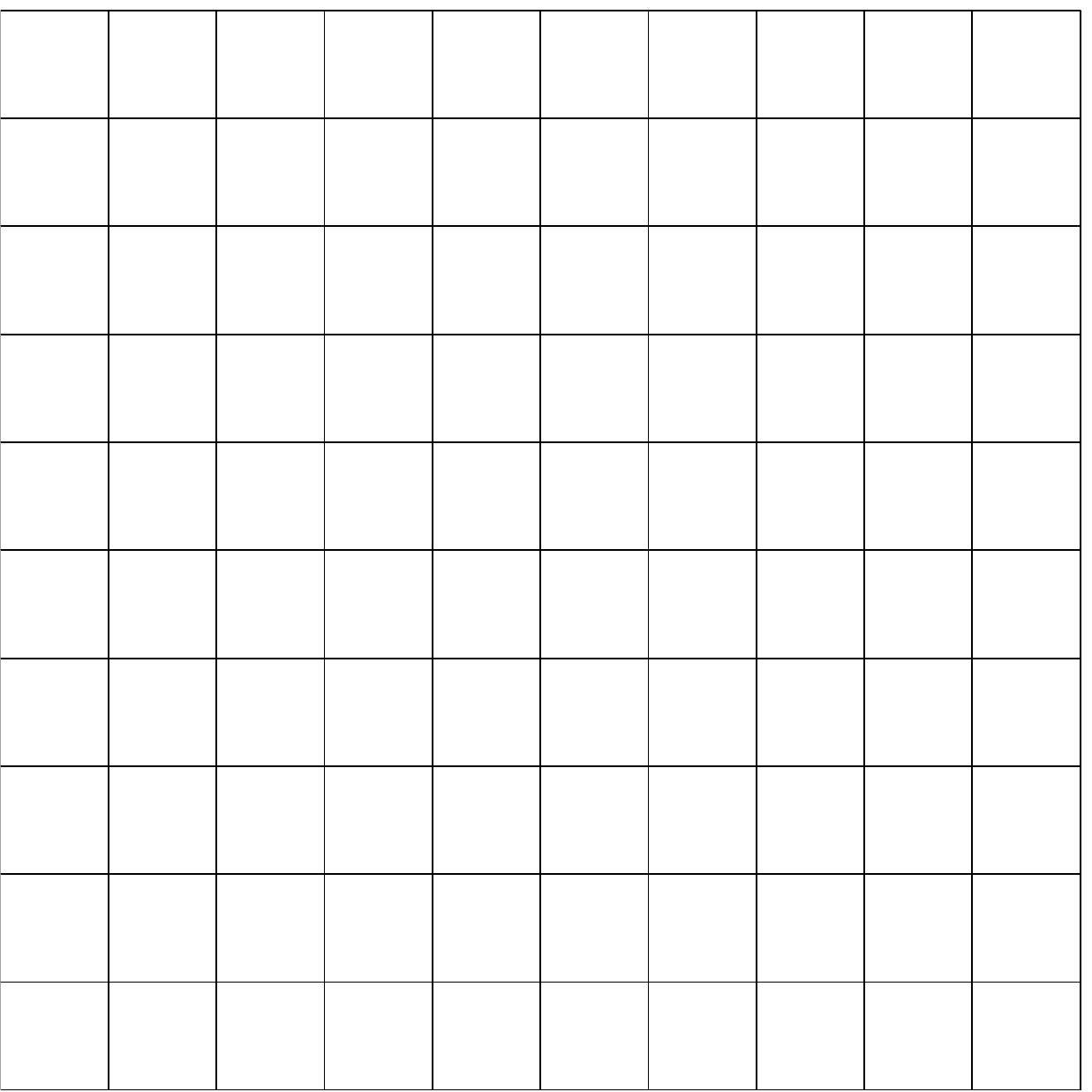}
  \end{center}
  \caption{\texttt{square} mesh}
 \label{fig:square}
 \end{minipage}
\hfill
 \begin{minipage}[b]{0.45\textwidth}
  \begin{center}
  \includegraphics[width=\textwidth]{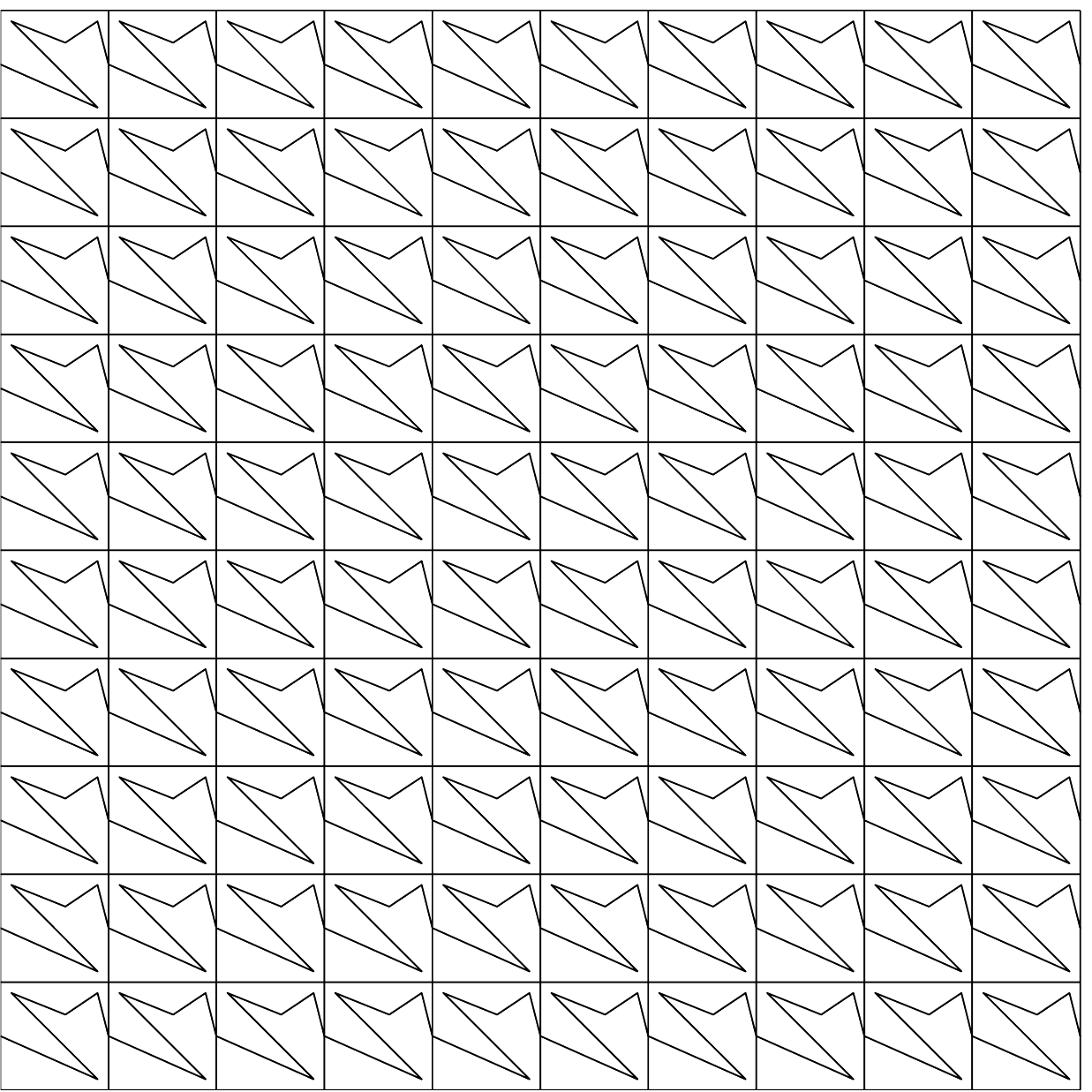}
  \end{center}
  \caption{\texttt{concave} mesh}
 \label{fig:concave}
 \end{minipage}
\hfill
\end{figure}

The third sequence of meshes (labelled \texttt{square}) is simply a decomposition
of the domain in 25, 100, 400 and 1600 equal squares, while the fourth sequence
(labelled \texttt{concave}) is obtained from the previous one by subdividing
each small square into two non-convex (quite nasty) polygons.
As before, the second meshes of the two sequences are shown in Fig. \ref{fig:square}
and in Fig. \ref{fig:concave} respectively.

\subsection{Convergence curves}

\begin{figure}
\hfill
 \begin{minipage}[b]{0.49\textwidth}
  \begin{center}
  \includegraphics[width=\textwidth]{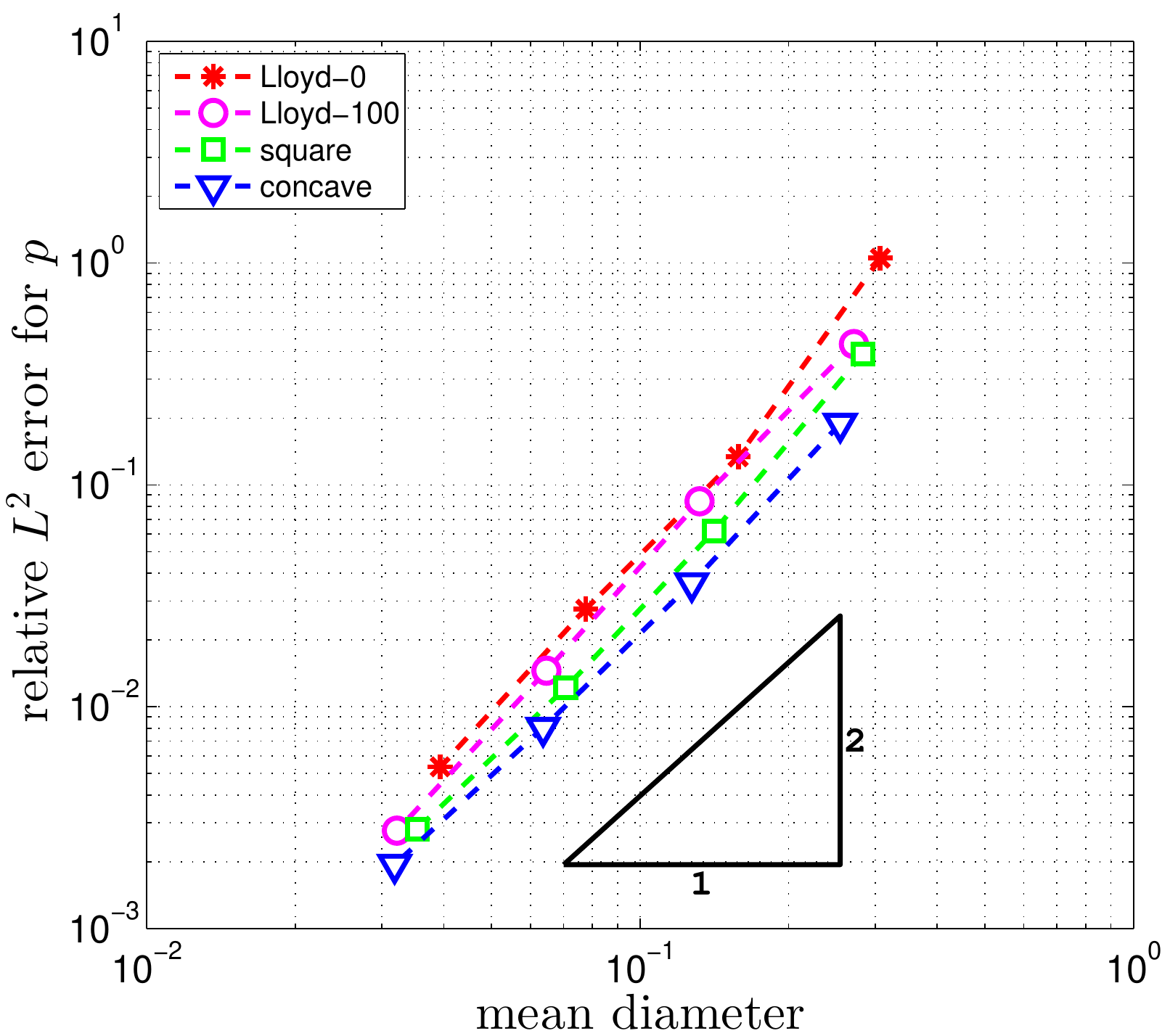}
  \end{center}
  \caption{$k=1$, relative $L^2$ error for $\ph$}
 \label{fig:err-p-k=1}
 \end{minipage}
\hfill
 \begin{minipage}[b]{0.49\textwidth}
  \begin{center}
  \includegraphics[width=\textwidth]{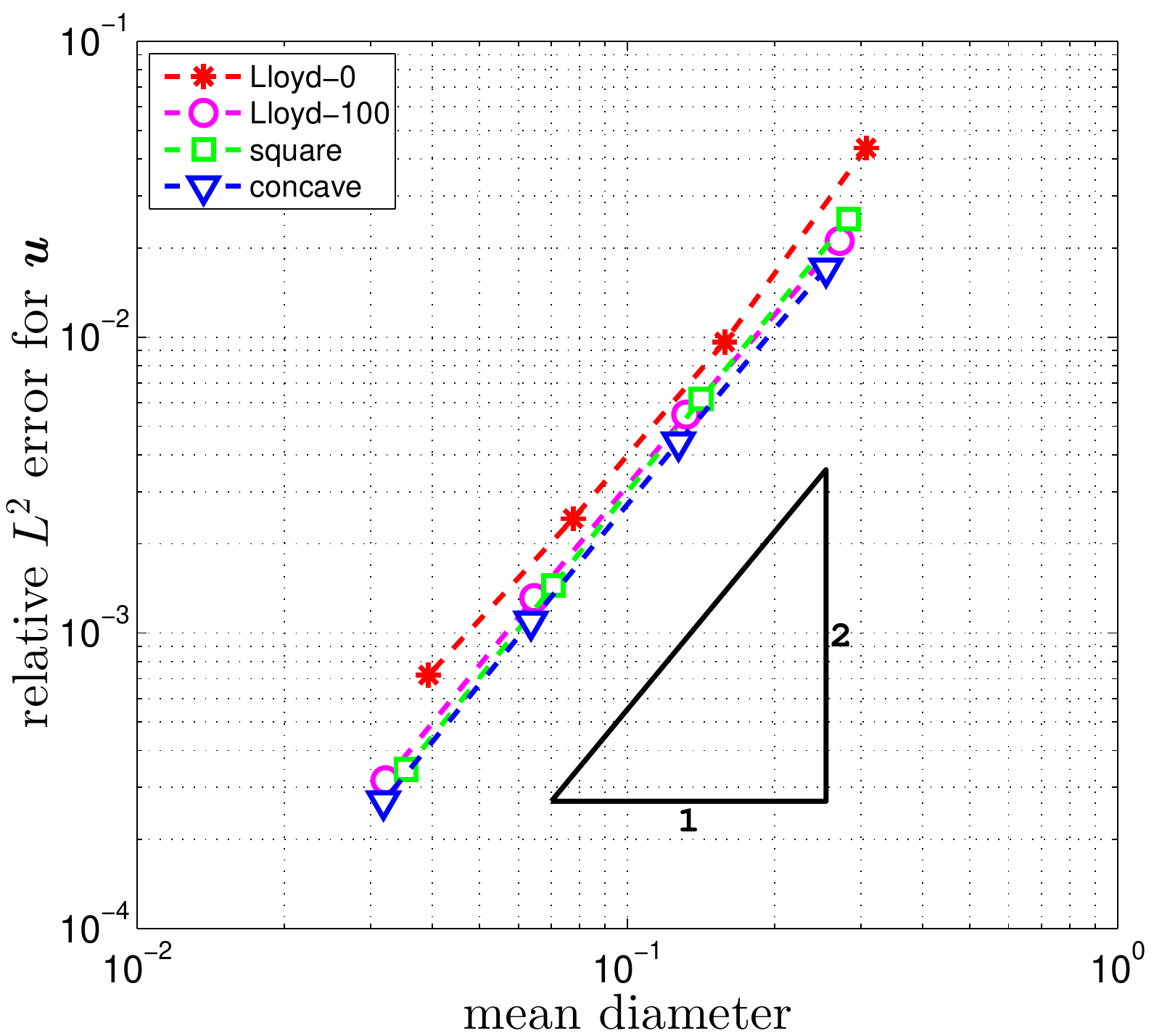}
  \end{center}
  \caption{$k=1$, relative $L^2$ error for $\uh$}
 \label{fig:err-u-k=1}
 \end{minipage}
\hfill
\end{figure}

\begin{figure}
\hfill
 \begin{minipage}[b]{0.49\textwidth}
  \begin{center}
  \includegraphics[width=\textwidth]{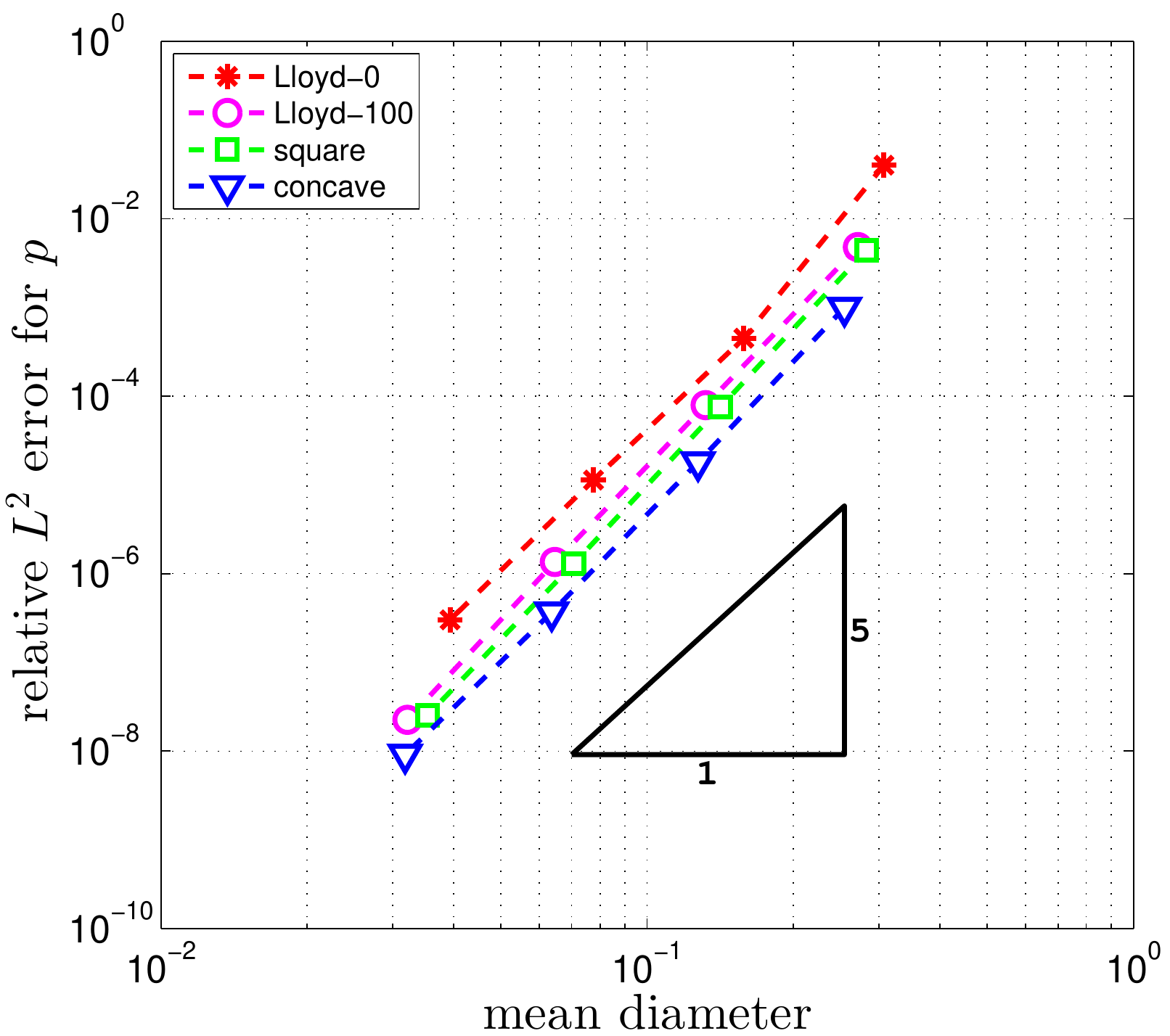}
  \end{center}
  \caption{$k=4$, relative $L^2$ error for $\ph$}
 \label{fig:err-p-k=4}
 \end{minipage}
\hfill
 \begin{minipage}[b]{0.49\textwidth}
  \begin{center}
  \includegraphics[width=\textwidth]{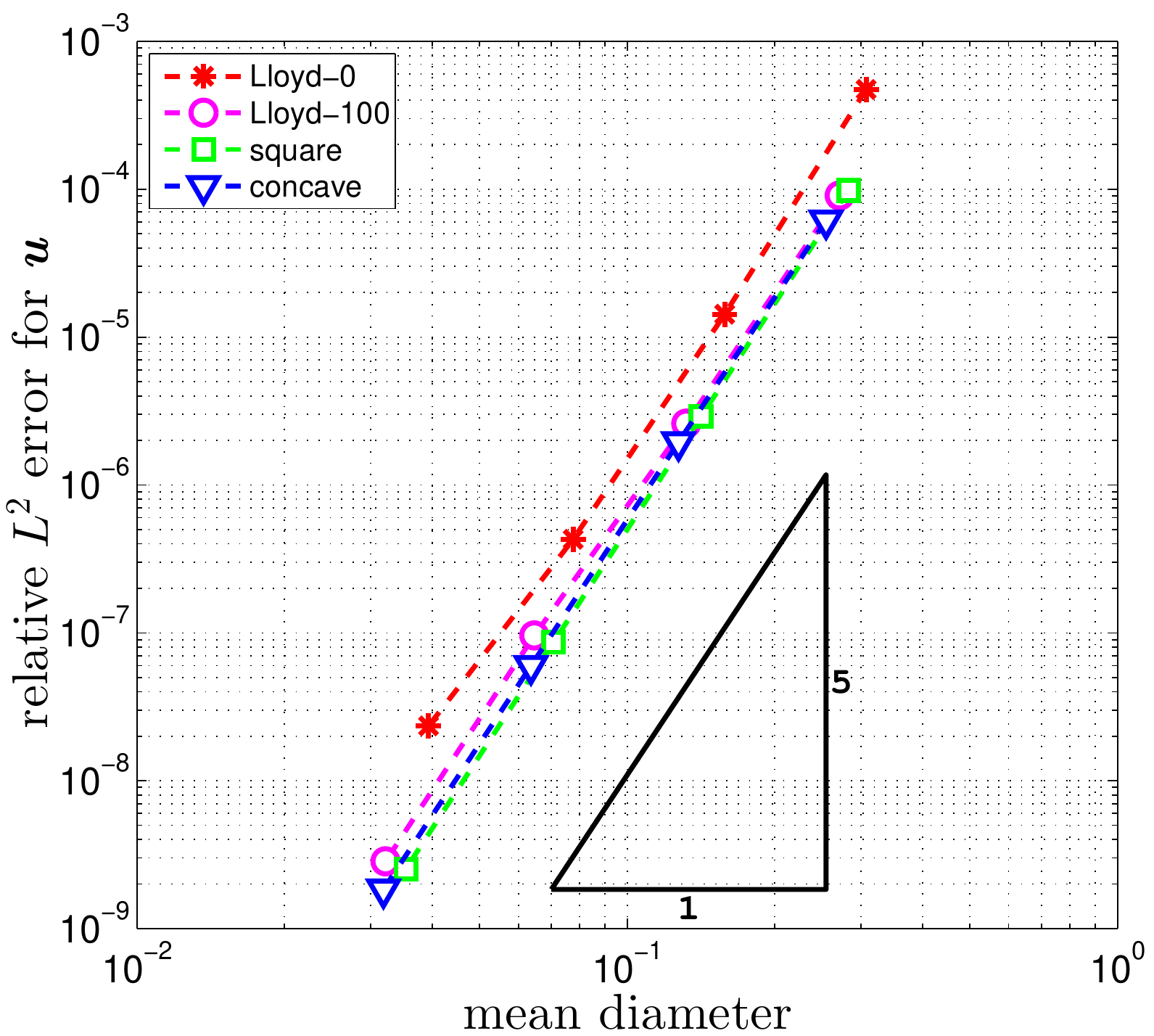}
  \end{center}
  \caption{$k=4$, relative $L^2$ error for $\uh$}
 \label{fig:err-u-k=4}
 \end{minipage}
\hfill
\end{figure}

%We start to show the convergence results for $k=1$.
In Figs. \ref{fig:err-p-k=1} and \ref{fig:err-u-k=1}
we report the relative error in $L^2$ for $\ph$ and $\uh$ respectively,
for the four mesh sequences in the case $k=1$.
In Figs.~\ref{fig:err-p-k=4} and \ref{fig:err-u-k=4}
we show the same convergence results for $k=4$.
%
%\begin{remark}

\begin{figure}
\hfill
 \begin{minipage}[b]{0.49\textwidth}
  \begin{center}
  \includegraphics[width=\textwidth]{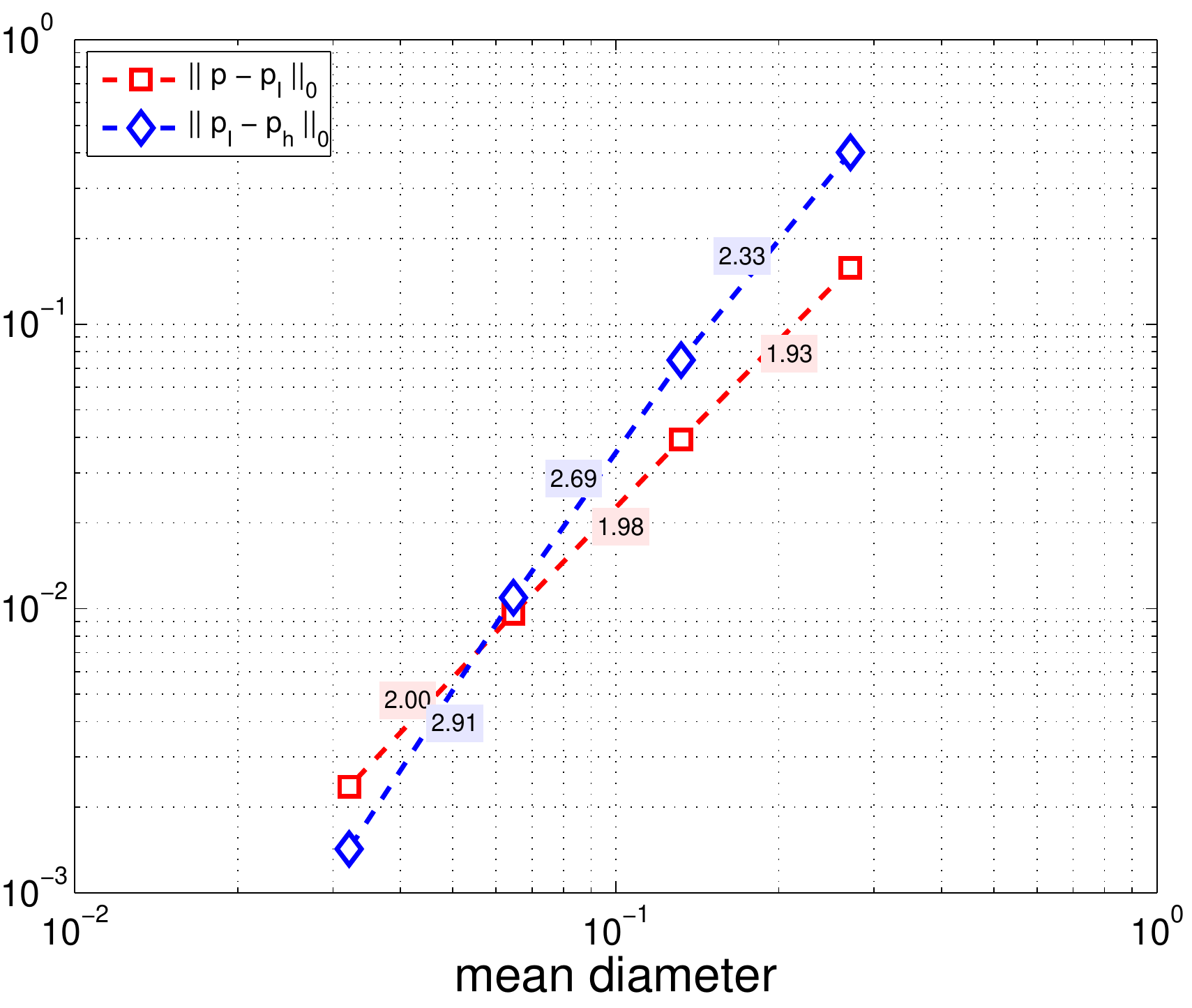}
  \end{center}
  \caption{$k=1$, superconvergence}
 \label{fig:super-k=1}
 \end{minipage}
\hfill
 \begin{minipage}[b]{0.49\textwidth}
  \begin{center}
  \includegraphics[width=\textwidth]{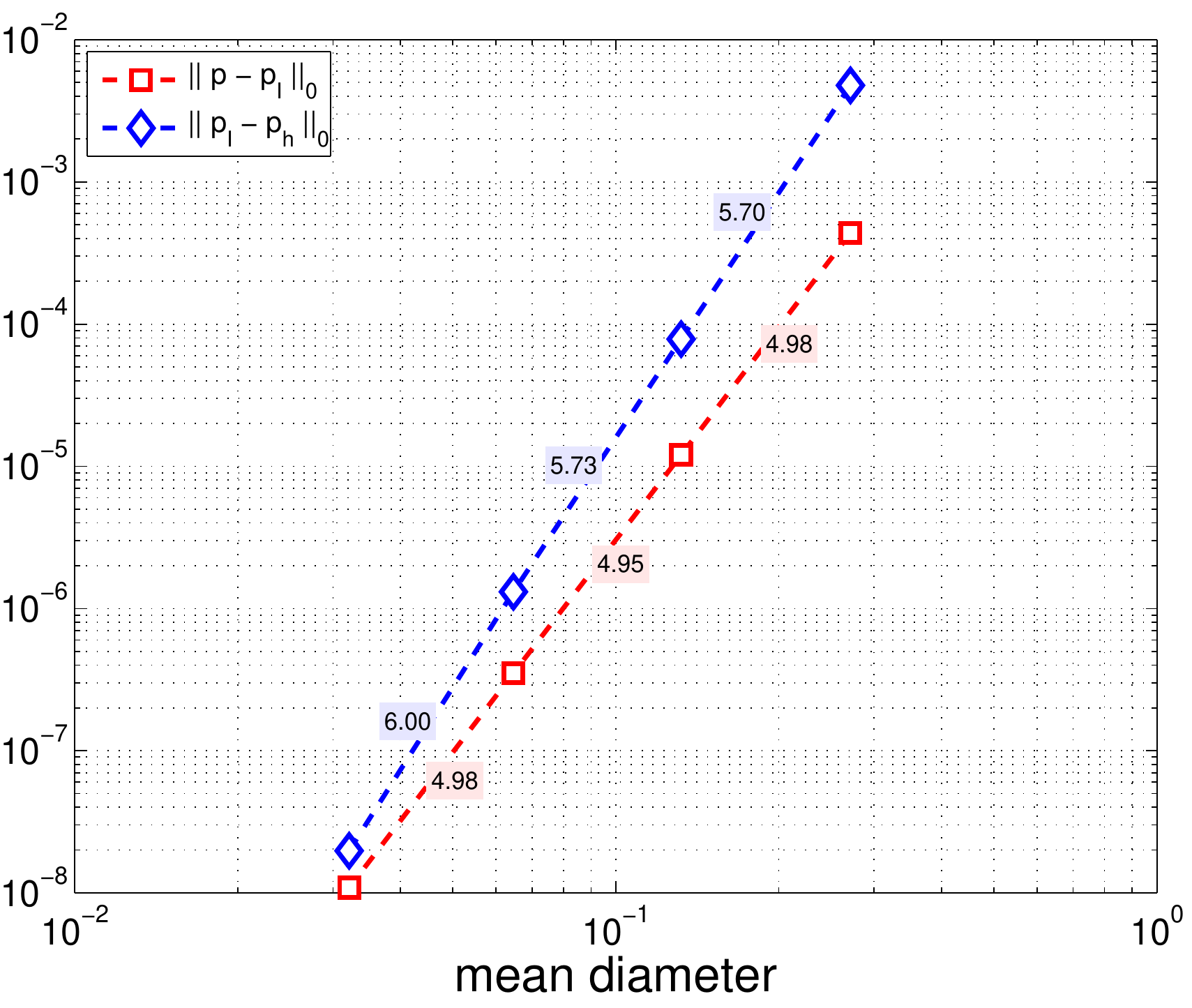}
  \end{center}
  \caption{$k=4$, superconvergence}
 \label{fig:super-k=4}
 \end{minipage}
\hfill
\end{figure}

A closer inspection of the convergence curves for the $L^2$ error between $\pex$ and $\ph$
shown in Figs.~\ref{fig:err-p-k=1} and \ref{fig:err-p-k=4} reveals that the slope is slightly
larger than expected for the coarsest meshes.
This behavior can be explained in following way. The $L^2$ error $\|p-\ph\|_{0}$
can be written as
\begin{equation}
\|p-\ph\|_{0}^2 = \|p-p_I\|_{0}^2 + \|p_I-\ph\|_{0}^2
\end{equation}
where we recall that on each element $p_I=\Pi^0_k p$.
As shown in Section \ref{superconvergence}, there is a superconvergence of $\ph$ to $p_I$:
\begin{equation}
\|p_I-\ph\|_{0}\leq C h^{k+2}.
\end{equation}
Hence, as long as $\|p_I-\ph\|_{0}$ is the dominant term in the error, we observe
a slope of $k+2$; when $h$ becomes smaller, the term $\|p-p_I\|_{0}$ takes over
and the slope becomes $k+1$ as expected. This is clearly shown in
Figs.~\ref{fig:super-k=1} and \ref{fig:super-k=4}
where $p-p_I$ and $\p_I-\ph$ are plotted in the case
of the \texttt{lloyd-100} meshes with $k=1$ and $k=4$, respectively.
For the sake of clarity, on each curve we have reported its slope.
%\end{remark}
%

We conclude that the Virtual Element Method behaves as expected and shows a remarkable
stability with respect to the shape of the mesh polygons.

\bibliographystyle{amsplain}

\bibliography{advdiff-bibliography}

\providecommand{\bysame}{\leavevmode\hbox to3em{\hrulefill}\thinspace}
\providecommand{\MR}{\relax\ifhmode\unskip\space\fi MR }
% \MRhref is called by the amsart/book/proc definition of \MR.
\providecommand{\MRhref}[2]{%
  \href{http://www.ams.org/mathscinet-getitem?mr=#1}{#2}
}
\providecommand{\href}[2]{#2}
\begin{thebibliography}{10}

\bibitem{projectors}
B.~Ahmad, A.~Alsaedi, F.~Brezzi, L.~D. Marini, and A.~Russo, \emph{Equivalent
  projectors for virtual element methods}, Comput. Math. Appl. \textbf{66}
  (2013), no.~3, 376--391.

\bibitem{ABMV14}
P.~F. Antonietti, L.~Beir{\~a}o~da Veiga, D.~Mora, and M.~Verani, \emph{A
  stream virtual element formulation of the {S}tokes problem on polygonal
  meshes}, SIAM J. Numer. Anal. \textbf{52} (2014), no.~1, 386--404.

\bibitem{Antonietti}
P.~F. Antonietti, N.~Bigoni, and M.~Verani, \emph{Mimetic discretizations of
  elliptic control problems}, J. Sci. Comput. \textbf{56} (2013), no.~1,
  14--27.

\bibitem{ABCM}
D.~N. Arnold, F.~Brezzi, B.~Cockburn, and L.~D. Marini, \emph{Unified analysis
  of discontinuous {G}alerkin methods for elliptic problems}, SIAM J. Numer.
  Anal. \textbf{39} (2001), no.~5, 1749--1779.

\bibitem{Arroyo-Ortiz}
M.~Arroyo and M.~Ortiz, \emph{Local maximum-entropy approximation schemes},
  Meshfree methods for partial differential equations {III}, Lect. Notes
  Comput. Sci. Eng., vol.~57, Springer, Berlin, 2007, pp.~1--16.

\bibitem{GFEM13}
I.~Babu\v{s}ka, U.~Banerjee, and J.~E. Osborn, \emph{Generalized finite element
  methods -- main ideas, results and perspective}, Int. J. Comput. Methods
  \textbf{01} (2004), no.~01, 67--103.

\bibitem{MESHLESS16}
I.~Babu\v{s}ka and J.~M. Melenk, \emph{The partition of unity method},
  Internat. J. Numer. Methods Engrg. \textbf{40} (1997), no.~4, 727--758.

\bibitem{GFEM17}
I.~Babu\v{s}ka and J.~E. Osborn, \emph{Generalized finite element methods:
  their performance and their relation to mixed methods}, SIAM J. Numer. Anal.
  \textbf{20} (1983), no.~3, 510--536.

\bibitem{Beirao-apos}
L.~Beir{\~a}o~da Veiga, \emph{A residual based error estimator for the mimetic
  finite difference method}, Numer. Math. \textbf{108} (2008), 387--406.

\bibitem{volley}
L.~Beir{\~a}o~da Veiga, F.~Brezzi, A.~Cangiani, G.~Manzini, L.~D. Marini, and
  A.~Russo, \emph{Basic principles of virtual element methods}, Math. Models
  Methods Appl. Sci. \textbf{23} (2013), no.~1, 199--214.

\bibitem{VEM-elasticity}
L.~Beir{\~a}o~da Veiga, F.~Brezzi, and L.~D. Marini, \emph{Virtual elements for
  linear elasticity problems}, SIAM J. Numer. Anal. \textbf{51} (2013), no.~2,
  794--812.

\bibitem{super-misti}
L.~Beir{\~a}o~da Veiga, F.~Brezzi, L.~D. Marini, and A.~Russo,
  \emph{${H}(\div)$ and ${H}(\curl)$-conforming {VEM}}, submitted.

\bibitem{variable-primal}
\bysame, \emph{Virtual element methods for general second order elliptic
  problems}, submitted.

\bibitem{hitchhikers}
\bysame, \emph{The hitchhiker's guide to the virtual element method}, Math.
  Models Methods Appl. Sci. \textbf{24} (2014), no.~8, 1541--1573.

\bibitem{MFD7}
L.~Beir{\~a}o~da Veiga, V.~Gyrya, K.~Lipnikov, and G.~Manzini, \emph{Mimetic
  finite difference method for the {S}tokes problem on polygonal meshes}, J.
  Comput. Phys. \textbf{228} (2009), no.~19, 7215--7232.

\bibitem{BLM09}
L.~Beir{\~a}o~da Veiga, K.~Lipnikov, and G.~Manzini, \emph{Convergence analysis
  of the high-order mimetic finite difference method}, Numer. Math.
  \textbf{113} (2009), no.~3, 325--356.

\bibitem{BLM11}
\bysame, \emph{Arbitrary-order nodal mimetic discretizations of elliptic
  problems on polygonal meshes}, SIAM J. Numer. Anal. \textbf{49} (2011),
  no.~5, 1737--1760.

\bibitem{MFD23}
\bysame, \emph{The mimetic finite difference method for elliptic problems},
  MS\&A. Modeling, Simulation and Applications, vol.~11, Springer-Verlag, 2014.

\bibitem{BeiraodaVeiga:Manzini:2008b}
L.~Beir{\~a}o~da Veiga and G.~Manzini, \emph{A higher-order formulation of the
  mimetic finite difference method}, SIAM J. Sci. Comput. \textbf{31} (2008),
  no.~1, 732--760.

\bibitem{BM13}
\bysame, \emph{A virtual element method with arbitrary regularity}, IMA J.
  Numer. Anal. \textbf{34} (2014), no.~2, 759--781.

\bibitem{Berrone-VEM}
M.~F. Benedetto, S.~Berrone, S.~Pieraccini, and S.~Scial{\`o}, \emph{The
  virtual element method for discrete fracture network simulations}, Comput.
  Methods Appl. Mech. Engrg. \textbf{280} (2014), 135--156.

\bibitem{Wriggers-2014}
S.~Biabanaki, A.~Khoei, and P.~Wriggers, \emph{Polygonal finite element methods
  for contact-impact problems on non-conformal meshes}, In press on {CMAME},
  DOI:10.1016/j.cma.2013.10.025, 198-221.

\bibitem{Bishop}
J.~E. Bishop, \emph{A displacement-based finite element formulation for general
  polyhedra using harmonic shape functions}, Internat. J. Numer. Methods Engrg.
  \textbf{97} (2014), no.~1, 1--31.

\bibitem{Bochev-Hyman}
P.~B. Bochev and J.~M. Hyman, \emph{Principles of mimetic discretizations of
  differential operators}, Compatible spatial discretizations, IMA Vol. Math.
  Appl., vol. 142, Springer, New York, 2006, pp.~89--119.

\bibitem{Bonelle:Ern:2014}
J.~Bonelle and A.~Ern, \emph{Analysis of compatible discrete operator schemes
  for elliptic problems on polyhedral meshes}, ESAIM Math. Model. Numer. Anal.
  \textbf{48} (2014), no.~2, 553--581.

\bibitem{Brezzi:Buffa:Lipnikov:2009}
F.~Brezzi, A.~Buffa, and K.~Lipnikov, \emph{Mimetic finite differences for
  elliptic problems}, M2AN Math. Model. Numer. Anal. \textbf{43} (2009), no.~2,
  277--295.

\bibitem{Maxwell-MFD}
F.~Brezzi, A.~Buffa, K.~Lipnikov, and G.~Manzini, \emph{The mimetic finite
  difference method for the 3d magnetostatic field problems on polyhedral
  meshes}, J. Comput. Phys. \textbf{230} (2011), 305--328.

\bibitem{BFM}
F.~Brezzi, R.~S. Falk, and L.~D. Marini, \emph{Basic principles of mixed
  virtual element methods}, ESAIM Math. Model. Numer. Anal. \textbf{48} (2014),
  no.~4, 1227--1240.

\bibitem{Brezzi:Lipnikov:Shashkov:2005}
F.~Brezzi, K.~Lipnikov, and M.~Shashkov, \emph{Convergence of mimetic finite
  difference method for diffusion problems on polyhedral meshes with curved
  faces}, Math. Models Methods Appl. Sci. \textbf{16} (2006), no.~2, 275--297.

\bibitem{Brezzi:Lipnikov:Shashkov:2006}
\bysame, \emph{Convergence of mimetic finite difference method for diffusion
  problems on polyhedral meshes with curved faces}, Math. Models Methods Appl.
  Sci. \textbf{16} (2006), no.~2, 275--297.

\bibitem{Brezzi:Lipnikov:Shashkov:Simoncini:2007}
F.~Brezzi, K.~Lipnikov, M.~Shashkov, and V.~Simoncini, \emph{A new
  discretization methodology for diffusion problems on generalized polyhedral
  meshes}, Comput. Methods Appl. Mech. Engrg. \textbf{196} (2007), no.~37-40,
  3682--3692.

\bibitem{Brezzi:Lipnikov:Simoncini:2005}
F.~Brezzi, K.~Lipnikov, and V.~Simoncini, \emph{A family of mimetic finite
  difference methods on polygonal and polyhedral meshes}, Math. Models Methods
  Appl. Sci. \textbf{15} (2005), no.~10, 1533--1551.

\bibitem{Brezzi:Marini:plates}
F.~Brezzi and L.~D. Marini, \emph{Virtual element methods for plate bending
  problems}, Comput. Methods Appl. Mech. Engrg. \textbf{253} (2013), 455--462.

\bibitem{hourglass}
A.~Cangiani, G.~Manzini, A.~Russo, and N.~Sukumar, \emph{Hourglass
  stabilization and the virtual element method}, to appear, 2015.

\bibitem{chessa-S-B}
J.~Chessa, P.~Smolinski, and T.~Belytschko, \emph{The extended finite element
  method (xfem) for solidification problems}, Internat. J. Numer. Methods
  Engrg. \textbf{53} (2002), 1959--1977.

\bibitem{Paulino-nonlinear-polygonal}
H.~Chi, C.~Talischi, O.~Lopez-Pamies, and G.H. Paulino, \emph{Polygonal finite
  elements for finite elasticity}, In press on {INJME}. DOI: 10.1002/nme.4802.

\bibitem{Ciarlet-78}
P.G. Ciarlet, \emph{The finite element method for elliptic problems}, Studies
  in Mathematics and its Applications, vol.~4, North-Holland Publishing Co.,
  Amsterdam-New York-Oxford, 1978, 1978.

\bibitem{Cockburn-IMU}
B.~Cockburn, \emph{The hybridizable discontinuous {G}alerkin methods},
  Proceedings of the {I}nternational {C}ongress of {M}athematicians. {V}olume
  {IV}, Hindustan Book Agency, New Delhi, 2010, pp.~2749--2775.

\bibitem{Cockburn-Jay-Lazarov}
B.~Cockburn, J.~Gopalakrishnan, and R.~Lazarov, \emph{Unified hybridization of
  discontinuous {G}alerkin, mixed, and continuous {G}alerkin methods for second
  order elliptic problems}, SIAM J. Numer. Anal. \textbf{47} (2009), no.~2,
  1319--1365.

\bibitem{Cockburn-Jay-Sayas}
B.~Cockburn, J.~Gopalakrishnan, and F.-J. Sayas, \emph{A projection-based error
  analysis of {HDG} methods}, Math. Comp. \textbf{79} (2010), no.~271,
  1351--1367.

\bibitem{Cockburn-Guzman-HWang}
B.~Cockburn, J.~Guzm{\'a}n, and H.~Wang, \emph{Superconvergent discontinuous
  {G}alerkin methods for second-order elliptic problems}, Math. Comp.
  \textbf{78} (2009), no.~265, 1--24.

\bibitem{DiPietro-Ern-1}
D.~Di~Pietro and A.~Alexandre~Ern, \emph{A hybrid high-order locking-free
  method for linear elasticity on general meshes}, Comput. Methods Appl. Mech.
  Engrg. \textbf{283} (2015), no.~0, 1--21.

\bibitem{DiPietro-Ern}
D.~Di~Pietro and A.~Ern, \emph{Mathematical aspects of discontinuous {G}alerkin
  methods}, Math\'ematiques \& Applications (Berlin) [Mathematics \&
  Applications], vol.~69, Springer, Heidelberg, 2012.

\bibitem{DiPietro-Ern-2}
\bysame, \emph{{A family of arbitrary-order mixed methods for heterogeneous
  anisotropic diffusion on general meshes}},
  https://hal.archives-ouvertes.fr/hal-00918482, December 2013.

\bibitem{DiPietro-Ern-3}
\bysame, \emph{Hybrid high-order methods for variable-diffusion problems on
  general meshes}, in press, 2014.

\bibitem{DiPietro-Ern--Lemaire}
D.~Di~Pietro, A.~Ern, and S.~Lemaire, \emph{An arbitrary-order and
  compact-stencil discretization of diffusion on general meshes based on local
  reconstruction operators}, Comput. Methods Appl. Math. \textbf{14} (2014),
  no.~4, 461--472.

\bibitem{DR-85}
J.~Douglas, Jr. and J.~E. Roberts, \emph{Global estimates for mixed methods for
  second order elliptic equations}, Math. Comp. \textbf{44} (1985), no.~169,
  39--52. \MR{771029 (86b:65122)}

\bibitem{Droniou:Eymard:Gallouet:Herbin:2010}
J.~Droniou, R.~Eymard, T.~Gallou{\"e}t, and R.~Herbin, \emph{A unified approach
  to mimetic finite difference, hybrid finite volume and mixed finite volume
  methods}, Math. Models Methods Appl. Sci. \textbf{20} (2010), no.~2,
  265--295.

\bibitem{Droniou-gradient}
\bysame, \emph{Gradient schemes: a generic framework for the discretisation of
  linear, nonlinear and nonlocal elliptic and parabolic equations}, Math.
  Models Methods Appl. Sci. \textbf{23} (2013), no.~13, 2395--2432.

\bibitem{Droniou-SpIss}
Jerome Droniou, \emph{Finite volume schemes for diffusion equations:
  introduction to and review of modern methods}, Math. Models Methods Appl.
  Sci. \textbf{24} (2014), no.~8, 1575--1619. \MR{3200243}

\bibitem{Du:Faber99}
Q.~Du, V.~Faber, and M.~Gunzburger, \emph{Centroidal {V}oronoi tessellations:
  applications and algorithms}, SIAM Rev. \textbf{41} (1999), no.~4, 637--676.

\bibitem{Gillette-2}
M.~Floater, A.~Gillette, and N.~Sukumar, \emph{Gradient bounds for {W}achspress
  coordinates on polytopes}, SIAM J. Numer. Anal. \textbf{52} (2014), no.~1,
  515--532.

\bibitem{FHK06}
M.~Floater, K.~Hormann, and G.~K{\'o}s, \emph{A general construction of
  barycentric coordinates over convex polygons}, Advances in Computational
  Mathematics \textbf{24} (2006), no.~1-4, 311--331.

\bibitem{Floater-K-R}
M.~S. Floater, G.~K{\'o}s, and M.~Reimers, \emph{Mean value coordinates in 3d},
  Comput. Aided Geom. Design \textbf{22} (2005), 623--631.

\bibitem{Fries:Belytschko:2010}
T.-P. Fries and T.~Belytschko, \emph{The extended/generalized finite element
  method: an overview of the method and its applications}, Internat. J. Numer.
  Methods Engrg. \textbf{84} (2010), no.~3, 253--304.

\bibitem{Gain-PhD}
A.~L. Gain, \emph{Polytope-based topology optimization using a mimetic-inspired
  method}, Ph.D. thesis, University of Illinois at Urbana-Champaign, 2013.

\bibitem{Paulino-VEM}
A.~L. Gain, C.~Talischi, and G.~H. Paulino, \emph{On the {V}irtual {E}lement
  {M}ethod for three-dimensional linear elasticity problems on arbitrary
  polyhedral meshes}, Comput. Methods Appl. Mech. Engrg. \textbf{282} (2014),
  132--160.

\bibitem{XFEM84}
A.~Gerstenberger and W.~A. Wall, \emph{An extended finite element
  method/{L}agrange multiplier based approach for fluid-structure interaction},
  Comput. Methods Appl. Mech. Engrg. \textbf{197} (2008), no.~19-20,
  1699--1714.

\bibitem{HF}
Kai Hormann and Michael~S. Floater, \emph{Mean value coordinates for arbitrary
  planar polygons}, ACM Trans. Graph. \textbf{25} (2006), no.~4, 1424--1441.

\bibitem{MESHLESS24}
S.~R. Idelsohn, E.~O{\~n}ate, N.~Calvo, and F.~Del~Pin, \emph{The meshless
  finite element method}, Internat. J. Numer. Methods Engrg. \textbf{58}
  (2003), no.~6, 893--912.

\bibitem{MFD22}
K.~Lipnikov, G.~Manzini, and M.~Shashkov, \emph{Mimetic finite difference
  method}, J. Comput. Phys. \textbf{257} (2014), no.~part B, 1163--1227.

\bibitem{VEM19}
G.~Manzini, A.~Russo, and N.~Sukumar, \emph{New perspectives on polygonal and
  polyhedral finite element methods}, Math. Models Methods Appl. Sci.
  \textbf{24} (2014), no.~8, 1665--1699.

\bibitem{Manzini:Russo:Sukumar}
\bysame, \emph{New perspectives on polygonal and polyhedral finite element
  methods}, Math. Models Methods Appl. Sci. \textbf{24} (2014), no.~8,
  1665--1699.

\bibitem{POLY28}
S.~Martin, P.~Kaufmann, M.~Botsch, M.~Wicke, and M.~Gross, \emph{Polyhedral
  finite elements using harmonic basis functions.}, Comput. Graph. Forum
  \textbf{27} (2008), no.~5, 1521--1529.

\bibitem{GFEM-117}
J.M. Melenk and I.~Babuska, \emph{The partition of unity finite element method:
  basic theory and applications}, Comp. Methods Appl. Mech. Engrg. \textbf{139}
  (1996), 289--314.

\bibitem{Merle-Dolbow}
R.~Merle and J.~Dolbow, \emph{Solving thermal and phase change problems with
  the extended finite element method}, Comput. Mech. \textbf{28} (2002),
  339--350.

\bibitem{LIBROXFEM-128}
S.~Mohammadi, \emph{Extended finite element method}, Blackwell Publishing Ltd,
  2008.

\bibitem{VemSteklov}
D.~Mora, G.~Rivera, and R.~Rodr{\'\i}guez, \emph{A virtual element method for
  the {S}teklov eigenvalue problem}, CI2MA Pre-Publicaci{\'o}n 2014-27, in
  press on {M}ath. {M}od. {M}eth. {A}ppl. {M}ath., 2014.

\bibitem{Mu-Wang-Wei-Ye-Zhao}
L.~Mu, J.~Wang, G.~Wei, X.~Ye, and S.~Zhao, \emph{Weak {G}alerkin methods for
  second order elliptic interface problems}, J. Comput. Phys. \textbf{250}
  (2013), 106--125.

\bibitem{Mu-Wang-Ye-2}
L.~Mu, J.~Wang, and X.~Ye, \emph{A weak {G}alerkin finite element method with
  polynomial reduction}, arXiv:1304.6481, 2013.

\bibitem{Nguyen-Peraire-Cockburn}
N.~C. Nguyen, J.~Peraire, and B.~Cockburn, \emph{An implicit high-order
  hybridizable discontinuous galerkin method for linear convection-diffusion
  equations}, J. Comput. Phys. \textbf{228} (2009), no.~9, 3232--3254.

\bibitem{Oswald}
J.~Oswald, R.~Gracie, R.~Khare, and T.~Belytschko, \emph{An extended finite
  element method for dislocations in complex geometries: Thin films and
  nanotubes}, Comp. Methods Appl. Mech. Engrg \textbf{198} (2009), 1872--1886.

\bibitem{GFEM146}
T.~Rabczuk, S.~Bordas, and G.~Zi, \emph{On three-dimensional modelling of crack
  growth using partition of unity methods}, Computers \& Structures \textbf{88}
  (2010), no.~23–24, 1391 -- 1411, Special Issue: Association of
  Computational Mechanics – United Kingdom.

\bibitem{Gillette-1}
A.~Rand, A.~Gillette, and C.~Bajaj, \emph{Interpolation error estimates for
  mean value coordinates over convex polygons}, Advances in Computational
  Mathematics \textbf{39} (2013), no.~2, 327--347.

\bibitem{BEM-weisser}
S~Rjasanow and S.~Weisser, \emph{Fem with trefftz trial functions on polyhedral
  elements}, J. of Comp. and Appl. Math. \textbf{263} (2014), 202--217.

\bibitem{SVC}
B.G. Smith, B.L.~Jr. Vaughan, and D.L. Chopp, \emph{The extended finite element
  method for boundary layer problems in biofilm growth}, Comm. App. Math. and
  Comp. Sci. \textbf{2} (2007), 35--56.

\bibitem{spiegel}
M.~Spiegel~{et al.}, \emph{Tetrahedral vs. polyhedral mesh size evaluation on
  flow velocity and wall shear stress for cerebral hemodynamic simulation},
  Comp. Meth. in Biomech. and Biomed. Engrng. \textbf{14} (2011), 9--22.

\bibitem{Suku04}
N.~Sukumar, \emph{Construction of polygonal interpolants: a maximum entropy
  approach}, Internat. J. Numer. Methods Engrg. \textbf{61} (2004), no.~12,
  2159--2181.

\bibitem{suku-CMB}
N.~Sukumar, D.L. Chopp, N.~M{\"o}es, and T.~Belytschko, \emph{Modeling holes
  and inclusions by level sets in the extended finite-element method}, Comp.
  Methods Appl. Mech. Engrg \textbf{190} (2001), 6183--6200.

\bibitem{Sukumar:Malsch:2006}
N.~Sukumar and E.~A. Malsch, \emph{Recent advances in the construction of
  polygonal finite element interpolants}, Arch. Comput. Methods Engrg.
  \textbf{13} (2006), no.~1, 129--163.

\bibitem{CRACKXFEM180}
N.~Sukumar, N.~M{\"o}es, B.~Moran, and T.~Belytschko, \emph{Extended finite
  element method for three-dimensional crack modelling}, Internat. J. Numer.
  Methods Engrg. \textbf{48} (2000), no.~11, 1549--1570.

\bibitem{ST04}
N.~Sukumar and A.~Tabarraei, \emph{Conforming polygonal finite elements},
  Internat. J. Numer. Methods Engrg. \textbf{61} (2004), no.~12, 2045--2066.

\bibitem{Sutradhar-Paulino-Miller-Nguyen}
A.~Sutradhar, G.~H. Paulino, M.~J. Miller, and T.~H. Nguyen, \emph{Topology
  optimization for designing patient-specific large craniofacial segmental bone
  replacements}, Proc. Natl. Acad. Sci. U.S.A \textbf{107} (2010),
  13222--13227.

\bibitem{POLY37}
C.~Talischi and G.~H. Paulino, \emph{Addressing integration error for polygonal
  finite elements through polynomial projections: a patch test connection},
  Math. Models Methods Appl. Sci. \textbf{24} (2014), no.~8, 1701--1727.

\bibitem{TPPM10}
C.~Talischi, G.~H. Paulino, A.~Pereira, and I.~F.~M. Menezes, \emph{Polygonal
  finite elements for topology optimization: A unifying paradigm}, Internat. J.
  Numer. Methods Engrg. \textbf{82} (2010), no.~6, 671--698.

\bibitem{vigneron}
L.M. Vigneron, J.G. Verly, and S.K. Warfield, \emph{On extended finite element
  method (xfem) for modelling of organ deformations associated with surgical
  cuts}, S. Cotin and D. Metaxas, editors, Medical Simulation, volume 3078 of
  Lecture Notes in Computer Science, Springer, Berlin, 2004.

\bibitem{Wachspress75}
E.~Wachspress, \emph{A rational finite element basis}, Academic Press, Inc.,
  New York-London, 1975, Mathematics in Science and Engineering, Vol. 114.

\bibitem{WMLB}
G.J. Wagner, N.~M{\"o}es, W.K. Liu, and T.~Belytschko, \emph{The extended
  finite element method for rigid particles in stokes flow}, Internat. J.
  Numer. Methods Engrg \textbf{51} (2991), 293--313.

\bibitem{Wang-1}
J.~Wang and X.~Ye, \emph{A weak {G}alerkin finite element method for
  second-order elliptic problems}, J. Comput. Appl. Math. \textbf{241} (2013),
  103--115.

\bibitem{Wang-2}
\bysame, \emph{A weak {G}alerkin mixed finite element method for second order
  elliptic problems}, Math. Comp. \textbf{83} (2014), no.~289, 2101--2126.

\bibitem{POLY47}
J.~Warren, \emph{Barycentric coordinates for convex polytopes}, Advances in
  Computational Mathematics \textbf{6} (1996), no.~1, 97--108.

\end{thebibliography}

\end{document}